%% file: Maximizing_measure_for_CS_via_Blur_shift_spaces.tex
\documentclass[cupthm,crop,info]{CUP-JNL-ETS}%

\usepackage{graphicx}
\usepackage{multicol,multirow}
\usepackage{amsmath,amssymb,amsfonts}
\usepackage{mathrsfs}
\usepackage{rotating}
\usepackage{appendix}
\usepackage[numbers]{natbib}
\usepackage{amsthm}
\usepackage{lipsum}

\theoremstyle{cupplain}
\newtheorem{theorem}{Theorem}[section]
\newtheorem{lemma}[theorem]{Lemma}
\newtheorem{corollary}[theorem]{Corollary}
\theoremstyle{cupdefinition}
\newtheorem{definition}{Definition}[section]
\theoremstyle{cupremark}
\newtheorem{remark}[theorem]{Remark}
\newtheorem{example}[theorem]{Example}
\theoremstyle{cupproof}
\newtheorem{proof}{Proof}

\numberwithin{equation}{section}

\input{packs.tex}
\jyear{2025}

\begin{document}

\begin{Frontmatter}

\title{Maximizing measures for countable alphabet shifts via blur shift spaces}

\author[1]{\gname{Eduardo} \sname{Garibaldi}}
\author[2]{\gname{Jo\~ao T. A.} \sname{Gomes}}
\author[3]{\gname{Marcelo} \sname{Sobottka}}

\address[1]{
\orgdiv{Institute of Mathematics, Statistics and Scientific Computing},
\orgname{University of Campinas},
\orgaddress{
\city{Campinas},
\postcode{13083-859},
\state{SP},
\country{Brazil}
}
\\
(\email{garibaldi@ime.unicamp.br})
}
\address[2]{
\orgdiv{Center of Exact and Technological Sciences},
\orgname{Federal University of Rec\^oncavo da Bahia},
\orgaddress{
\city{Cruz das Almas},
\postcode{44380-000},
\state{BA},
\country{Brazil}
} \\
(\email{jtagomes@ufrb.edu.br})
}
\address[3]{
\orgdiv{Department of Mathematics},
\orgname{Federal University of Santa Catarina},
\orgaddress{
\city{Florian\'opolis},
\postcode{88040-900},
\state{SC},
\country{Brazil}
}
\\
(\email{marcelo.sobottka@ufsc.br})
}


\maketitle

\authormark{E. Garibaldi, J. Gomes and M. Sobottka}

\titlemark{Maximizing measure for countable alphabet shifts via blur shift spaces}

\begin{abstract}
For upper semi-continuous potentials defined on shifts over countable alphabets, this paper ensures sufficient conditions for the existence of a maximizing measure.
We resort to the concept of blur shift, introduced by T.~Almeida and M.~Sobottka as a compactification method for countable alphabet shifts consisting of adding new symbols given by blurred subsets of the alphabet.
Our approach extends beyond the Markovian case to encompass more general countable alphabet shifts.
In particular, we guarantee a convex characterization and compactness for the set of blur invariant probabilities with respect to the discontinuous shift map.
\end{abstract}




\end{Frontmatter}

\section{Main Result}
\label{section-introduction}

Let $ \Sigma \subset \mathscr{A}^{\mathbb{N}} $ be any one-sided shift space defined over a countable alphabet $ \mathscr{A} $ and let $ \sigma \,:\, \Sigma \to \Sigma $ be the left shift map. 
To be concrete, we assume from this point forward that $ \mathscr{A} $ is identified with the set of non-negative integers $ \mathbb{Z}_{+} $.
As usual, $ \mathscr{A} $ is endowed with the discrete topology and $ \Sigma $ with the product topology.

We introduce the language $ \mathscr{L} $ of $ \Sigma $ as the set of allowed words, i.e., the union of the following sets
\begin{align*}
\mathscr{L}_{0} & = \left\{ \epsilon \right\}, \enspace \text{where } \epsilon \text{ is the empty word in } \mathscr{A}^{0}, \\
\text{and} \quad
\mathscr{L}_{k} & = \left\{ \left( x_0, \ldots, x_{k-1} \right) \in \mathscr{A}^k \,:\, x \in \Sigma \right\}, \enspace \text{for all } k \geq 1.
\end{align*}
As a notational convention, we use $ x $, $ y $ or $ z $ for (infinite) sequences
and $ u $, $ v $ or $ w$ for (finite) words from $ \mathscr{L} $.
In particular, the length of a word $ w \in \mathscr{L}_{k} $ is indicated as $ \ell (w) = k $.
We define the predecessor set as
\begin{equation*}
\mathscr{P} (w) = \{ a \in \mathscr{L}_{1} \,:\, a w  \in \mathscr{L} \},
\end{equation*}
and similarly we introduce the follower sets as
\begin{align*}
\mathscr{F}_1 \left( w \right) 
& = \{ b \in \mathscr{L}_{1} \,:\, w b \in \mathscr{L} \} 
\quad \text{and}
\\
\mathscr{F}_m \left( w \right) 
& = \{ b \in \mathscr{L}_{1} \,:\, w v b \in \mathscr{L} \text{ for some } v \in \mathscr{L}_{m-1} \},
\quad \text{for any } m \geq 2.
\end{align*}
For the empty word, we assume that $ \mathscr{P} (\epsilon) = \mathscr{F}_1 (\epsilon) = \mathscr{L}_1 $.

We call a potential any Borel function $ A \,:\, \Sigma \to \mathbb{R} $ which is bounded from above.
We introduce its ergodic maximizing constant as
\begin{equation*}
\beta ( A ) = \sup \left\{ \int_{\Sigma} \, A \, d \mu \,:\, \mu \mbox{ is a $ \sigma $-invariant probability measure} \right\}.
\end{equation*}
Obviously, $ \beta (A) \in \big(-\infty, \, \sup A \big] $. 
In these notes, we study the $ \sigma $-invariant probability measures $ \mu $ that maximize the integral of $ A $ over $ \Sigma $, i.e., $ \int A \, d \mu = \beta (A) $. 
Our main objective is to provide sufficient conditions (on the space and on the potential) for the existence of these measures, called maximizing probabilities.

\begin{theorem-E}  \hypertarget{EThm}{}
Let $ \Sigma $ be a shift over a countable alphabet that satisfies both
\begin{itemize}
\item \hypertarget{fcpa}{\itshape finite cyclic predecessor assumption}: \enspace 
$ \mathscr{P} (a) \cap \mathscr{F}_m (a) $ is finite for every $ a \in \mathscr{L}_1 $ and for all $ m \geq 1 $;
\item \hypertarget{dpm}{\itshape denseness of periodic measures}: \enspace
the set of ergodic probabilities supported on periodic orbits of $ \Sigma $ is (weak$^{\ast}$) dense among the $ \sigma $-invariant measures.
\end{itemize}
Then, every upper semi-continuous potential $ A $ fulfilling
\begin{equation*}
\limsup_{i \to \infty} \, \sup A|_{[i]} < \beta (A)
\end{equation*}
has a maximizing probability. 
Furthermore, if~$ \, \sup A \leq \beta (A) $, then there exists a common finite subshift that contains the support of any maximizing measure for~$ A $. 
\end{theorem-E}

To the best of our knowledge, this is an inaugural result for potentials defined on a non-compact shift that are less regular than summable variation ones, a ty\-pi\-cal scenario (see \cite{JMU:ETDS06, BG:BBMS10, BF:ETDS13}).
Moreover, the key assumption about the potentials is the requirement that, except occasionally in a finite subset of the alphabet, the supremum of their restrictions to cylinders of size $ 1 $ is strictly smaller than the respective maximizing ergodic constant.
In particular, this embraces the coercive framework and the conditions of oscillations considered in previous works \cite{JMU:ETDS06, BG:BBMS10, BF:ETDS13}.
With respect to the dynamics, the most comprehensive context addressed in the literature consists of topologically transitive countable alphabet Markov shifts.
Our \hyperlink{EThm}{Existence Theorem} also applies to more general countable alphabet shifts, taking into account situations beyond the Markovian case
(see Section~\ref{section-non-markov-example}).
Note that in the Markov setting, topological transitivity implies the \hyperlink{dpm}{denseness of periodic measures} 
(see, for instance, \cite{CS:IJM10}).
An equivalent formulation of the \hyperlink{fcpa}{finite cyclic predecessor assumption} was previously used in \cite{IV:JAM21} to study the set of invariant measures of topologically transitive countable Markov shifts.
Here this hypothesis has a central role in the proof of the (weak$^{\ast}$) compactness of the set of invariant probabilities for a compact extension $ \big( \hat{\Sigma}, \hat{\sigma} \big) $ of $ \left( \Sigma, \sigma \right) $, known as the blur shift~\cite{AS:BSM21}, whose construction will be described in detail in the \hyperlink{sec-blur-shift}{next section}.

The rest of the paper is organized as follows.
After reviewing the main properties and results of blur shifts in Section~\ref{section-blur-shift}, we study their respective invariant probabilities in the \hyperlink{sec-invariant-probabilities}{third section}.
We not only give a convex characterization (Proposition~\ref{proposition-decomposition-blur-invariant-measures}), but we also show that 
the \hyperlink{fcpa}{finite cyclic predecessor assumption} and 
the \hyperlink{dpm}{denseness of periodic measures} are sufficient conditions for compactness of the set of blur invariant measures (Corollary~\ref{corollary-blur-measure-compact}).
In Section~\ref{section-potentials}, we describe a minimal upper semi-continuous extension of the potential on the blur shift space. 
At this point, we can clearly present in the \hyperlink{sec-ergodic-optimization}{fifth section} the ergodic optimization problem for the blur shift and relate it to the initial one for the original countable alphabet shift.
Section~\ref{section-existence-maximizing-measure} is devoted to the proof of the \hyperlink{EThm}{Existence Theorem}.
In the \hyperlink{sec-non-markov-example}{final section}, we exhibit a countable alphabet non-Markovian shift that satisfies all the hypotheses required by our main result, demonstrating the applicability of our framework beyond the Markovian setting.

\hypertarget{sec-blur-shift}{}
\section{Blur Shift Space}
\label{section-blur-shift}
The blur shift space was recently introduced in~\cite{AS:BSM21} as a compactification scheme for classical shift spaces.
It is a unified presentation and generalization of 
the Ott-Tomforde-Willis shift space~\cite{OTW:NYJM14} and
the Gon\c{c}alves-Royer ultragraph shift space~\cite{GR:IMRN19},
which were successfully used
to study and extend results on the equivalence between the conjugacy of a subclass of Markov shift spaces and the isomorphism of the $C^\ast$-algebras associated with their graphs.

To the best of our knowledge, this note is a pioneer application of blur shift compactification method to ergodic optimization theory. 

Proofs of all the results just mentioned along this section may be found in~\cite{OTW:NYJM14, AS:BSM21}.
The reader interested in concrete examples of shift compactification may also consult these references.

Let $ \Sigma $ be a shift space over a countable alphabet.
A set $ \mathscr{V} = \left\{ B_1, \ldots, B_s \right\} \subset 2^{ \mathscr{A} } $ is said to be a finite resolution of blurred sets if it verifies:
$ B_r $ is infinite for each $ 1 \leq r \leq s $;
$ B_i \cap B_j $ is finite for all $ 1 \leq i \neq j \leq s $;
and $ \mathscr{L}_{1} \setminus \bigcup_{r=1}^{s} B_r $ is finite. 
We define the blur shift space $ \hat{\Sigma} \subset \left( \mathscr{A} \sqcup \mathscr{V} \right)^{ \mathbb{N} } $
with resolution $ \mathscr{V} $ associated with $ \Sigma $ as 
\begin{align*}
& \hat{\Sigma} \,=\, \Sigma \, \sqcup \, \partial_{ B_1 } \Sigma \, \sqcup \, \cdots \, \sqcup \, \partial_{ B_s } \Sigma, \\
& \text{where} \quad
\partial_{ B_r } \Sigma = 
\left\{ 
\left( w, B_r, B_r, \ldots \right)  \,:\, 
w \in \mathscr{L} \text { and } \mathscr{F}_1 (w) \cap B_r \text{ is infinite}
\right\},
\quad \text{for } 1 \leq r \leq s.
\end{align*}
Note that each $ \partial_{ B_r } \Sigma $ is non-empty, since the empty word $ \epsilon \in \mathscr{L} $ fulfill $ \mathscr{F}_1 ( \epsilon ) \cap B_r = B_r $. 
Besides, we identify $ \left( \epsilon, B_r, B_r, \ldots \right) $ with $ \left( B_r, B_r, \ldots \right) $.

Based on the original shift language, we consider the following decomposition of the set $ \hat{\Sigma} \setminus \Sigma $:
\begin{align*}
\hat{\mathcal{L}}_{0} 
& = \left\{ \left( B_r, B_r, \ldots \right) \in \hat{\Sigma} \,:\, 1 \leq r \leq s \right\} 
\\
\text{and} \quad
\hat{\mathcal{L}}_{k}  
& = \left\{ \left( w, B_r, B_r, \ldots \right) \in \hat{\Sigma} \,:\, w \in \mathscr{L}_k \text{ and } 1 \leq r \leq s \right\},
\quad\text{for } k \geq 1.
\end{align*}
For notational convenience, we denote  
$ \partial \Sigma 
= \hat{\Sigma} \setminus \Sigma 
= \bigsqcup_{r = 1}^{s} \partial_{ B_r } \Sigma 
= \bigsqcup_{k \geq 0} \hat{\mathcal{L}}_{k} $.

Let us recall how we can interpret $ \hat{\Sigma} $ as compact space in which $ \Sigma $ can naturally be seen as its topological subspace.
First, consider on $ \mathscr{A} \sqcup \mathscr{V} $ the topology generated by the family formed of singletons of $ \mathscr{A} $ and of sets of the form $ B_r \setminus S $, where $ S $ is some finite subset of $ B_r $ and $ 1 \leq r \leq s $.
Thus the set $ \big( \mathscr{A} \sqcup \mathscr{V} \big)^{ \mathbb{N} } $ is equipped with the corresponding product topology.
We say that $ x $, $ y \in \big( \mathscr{A} \sqcup \mathscr{V} \big)^{ \mathbb{N} } $ are equivalent, and we write $ x \sim y $, if
\begin{itemize}
\item either $ x $, $ y \in \mathscr{A}^{ \mathbb{N} } $ and $ x = y $;
\item or $ x $, $ y \in \big( \mathscr{A} \sqcup \mathscr{V} \big)^{ \mathbb{N} } \setminus \mathscr{A}^{ \mathbb{N} } $, 
$ x = ( u , B_i, \ldots) $, 
$ y = ( v , B_j, \ldots) $ 
with
$ u = v \in \bigsqcup_{ k \geq 0} \mathscr{A}^k $
and
$ B_i = B_j \in \mathscr{V} $.
\end{itemize}
From this equivalence relation, we can provide the quotient topology on $ \big( \mathscr{A} \sqcup \mathscr{V} \big)^{ \mathbb{N} } / { \sim } $.
Hence, since $ \hat{\Sigma} $ can be seen as a subset of $ \big( \mathscr{A} \sqcup \mathscr{V} \big)^{ \mathbb{N} } / { \sim } $,
the blur shift space $ \hat{\Sigma } $  inherits this quotient topology as a subspace.
It can be shown that $ \hat{\Sigma} $ is therefore a compact Hausdorff metrizable space (see Theorem~3.24 and Corollary~3.18 of \cite{AS:BSM21}).
The topological structure on $ \hat{\Sigma} $ has a countable basis of generalized cylinder sets, namely:
\begin{align*}
Z [ w ] & = \left\{ x \in \hat{\Sigma} \,:\, 
x_i = w_i, \text{ for all } 0 \leq i < \ell(w) \right\} 
\qquad\text{and} \\
Z [ w B_r; S ] & = \left\{ x \in \hat{\Sigma} \,:\, 
x_i = w_i, \text{ for all } 0 \leq i < \ell(w), \ 
x_{ \ell (w) } = B_r \text{ or } x_{ \ell (w) } \in B_r \setminus S \right\},
\end{align*}
where $ w \in \mathscr{L} $ is any finite word, $ B_r \in \mathscr{V} $ is some blurred set and $ S \subset B_r $ is a finite set.
In particular, the usual cylinder set 
$ [ w ] = \left\{ x \in \Sigma \,:\, 
x_i = w_i, \text{ for all } 0 \leq i < \ell(w) \right\}
\subset \Sigma $ 
coincides with the set $ Z [ w ] \cap \Sigma $.
It follows from this fact that the topology of $ \Sigma $ is the subspace topology from $ \hat{\Sigma} $ on $ \Sigma $.  
Note also that
$ Z [ \epsilon ] = \hat{\Sigma} $ and
$ Z [ B_r; S ] = \left\{ x \in \hat{\Sigma} \,:\, x_0 = B_r \text{ or } x_0 \in B_r \setminus S \right\} $. 
The generalized cylinders sets are simultaneously closed and open 
(see Proposition~3.1 of~\cite{AS:BSM21}).
The neighborhood basis for a point is
\begin{align*}
\left\{ Z [ x_0, \ldots, x_{n} ] \,:\, n \geq 0 \right\} \quad &  \text{if } \left( x_0, x_1, \ldots \right) \in \Sigma, \\
\left\{ Z [ w B_r; S ] \,:\, S \text{ finite subset of } B_r \right\} \quad & \text{if } \left( w, B_r, B_r, \ldots \right) \in \partial \Sigma \text{ with } \ell(w) \geq 1, \\
\text{or} \qquad
\left\{ Z [ B_r; S ] \,:\, S \text{ finite subset of } B_r \right\} \quad &  \text{if } \left( B_r, B_r, \ldots \right) \in \hat{\mathcal{L}}_{0}
\end{align*}
(see Proposition~3.2 of~\cite{AS:BSM21}).
The original shift space $ \Sigma $ is a dense subset in~$ \hat{\Sigma} $.
Indeed, given $ \left( w, B_r, B_r, \ldots \right) \in \partial \Sigma $ and some neighborhood $ Z [ w B_r, S ] $, since $ \left( \mathscr{F}_1 (w) \cap B_r \right) \setminus S $ is infinite, let $ a \in B_r \setminus S $ such that $ w a \in \mathscr{L} $.
In particular, there exists $ y \in \Sigma $ fulfilling $ \left( y_0, \ldots, y_{\ell(w)} \right) = w a $
and, hence, $ y \in Z [ w B_r, S ] $.

The convergence on $ \hat{\Sigma} $ can be described as follows.
\begin{lemma}[{\cite[Corollary~2.17]{OTW:NYJM14}}]  \label{lemma-convegence-hat-sigma} \ 
\begin{enumerate}[label=(\roman*)] 
\item\label{item-convegence-sigma}
The sequence $ \{ x^n \} \subset \hat{\Sigma} $ converges to $ x \in \Sigma $
if, and only if,
for every positive integer $ M $, there exists an integer $ N > 0 $ such that $ n > N $ implies $x^{n}_{i} = x_i $ for all $ 1 \leq  i \leq M $;
\item\label{item-convegence-partial-sigma}
The sequence $ \{ x^n \} \subset \hat{\Sigma} $ converges to $ \left( w, B_r, B_r, \ldots \right) \in \partial \Sigma $
if, and only if,
for every finite subset $ S \subset B_r $, there exists an integer $ N > 0 $ such that $ n > N $ implies 
$ x^{n}_{i} = w_i $ for all $ 0 \leq i < \ell ( w ) $ and
$ x^n_{ \ell (w) } = B_r $ or $ x^n_{ \ell (w) } \in B_r \setminus S $.
\end{enumerate}
\end{lemma}

The blur shift space $ \hat{\Sigma} $ obviously admits as measurable structure the Borel sigma-algebra generated by the generalized cylinder sets, which contains the Borel sigma-algebra of $ \Sigma $.
Let $ \mathcal{M} \left( \Sigma \right) $, $ \mathcal{M} \big( \hat{\Sigma} \big) $ denote, respectively, the sets of Borel probability measures on $ \Sigma $ and on $ \hat{\Sigma} $.
Both sets are convex and each one is considered equipped with the respective weak$^{\ast}$ topology. 
The following result allows us to commit the abuse of notation $ \mathcal{M} \left( \Sigma \right) \subset \mathcal{M} \big( \hat{\Sigma} \big) $.

\begin{lemma}  \label{lemma-topological-inclusion}
$ \mathcal{M} \left( \Sigma \right) $ is in a continuous bijective correspondence with the subset of $ \mathcal{M} \big( \hat{\Sigma} \big) $ formed by those probabilities that gives mass $ 1 $ to $ \Sigma $.
\end{lemma}

\begin{proof}
For $ \mu \in \mathcal{M} \big( \Sigma \big) $, note that 
$ \dot{\mu} \left( \,\cdot\, \right) 
= \mu \left( \,\cdot\, \cap \Sigma \right) $
defines a Borel probability on $ \hat{\Sigma} $.
This remark allows us to consider the map
$ \mu \in \mathcal{M} \big( \Sigma \big) 
\, \mapsto \, 
\Phi \left( \mu \right) = \dot{\mu} 
\in \mathcal{M} \big( \hat{\Sigma} \big) $, 
which is injective due to standard arguments on extensions of measures
(see, for instance, Theorem~3.3 of \cite{Bil:JWS17}).
Clearly, any $ \hat{\mu} \in \mathcal{M} \big( \hat{\Sigma} \big) $ with $ \hat{\mu} \left( \Sigma \right) = 1 $ may be seen as the image of its restriction on the Borel measurable sets of $ \Sigma $.
By definition and linearity, we have that
$ \int_{ \Sigma } \, \hat{\varphi}|_{ \Sigma } \, d \mu 
= \int_{ \hat{\Sigma} } \, \hat{\varphi} \, d \dot{\mu} $
for every simple function $ \hat{\varphi} $ on $ \hat{\Sigma} $,
and thus by the monotone convergence theorem, for every continuous function $ \hat{\varphi} $ on $ \hat{\Sigma} $.
This fact shows that if $ \hat{V} $ is a basis set of the weak$^{\ast}$ topology on $ \mathcal{M} \big( \hat{\Sigma} \big) $, then
$ \Phi^{-1} \big( \hat{V} \big) $ is a basis set of 
weak$^{\ast}$ topology on $ \mathcal{M} \left( \Sigma \right) $.
In other words, $ \Phi $ is continuous.
\end{proof}

As $ \hat{\Sigma} $ is compact metrizable space, it is well known that $ \mathcal{M} \big( \hat{\Sigma} \big) $ is a weak$^{\ast}$ compact metrizable space (see, for instance, Proposition~6.5 of \cite{Wal:SPR81}).
As a subset of $ \mathcal{M} \big( \hat{\Sigma} \big) $,  $ \mathcal{M} \big( \Sigma \big) $ is clearly not closed: for example, 
$ \delta_{ \left( k, 0, 0, \ldots \right) } 
\, \stackrel{\ast}{\rightharpoonup} \,
\delta_{ \left( B_r, B_r, \ldots \right) } $
as $ k \in B_r $ tends to~$ \infty $.

We now introduce the (blur) shift map $ \hat{\sigma} \,:\, \hat{\Sigma} \longrightarrow \hat{\Sigma} $ as the usual left shift map.
It~is obvious that $ \hat{\sigma} $ coincides with the shift map $ \sigma $ on $ \Sigma $.
Note that, in $ \partial_{  B_r } \Sigma $, there is a single fixed point 
$ \left( B_r, B_r, \ldots \right) = \hat{\sigma} \left( B_r, B_r, \ldots \right) $
which absorbs all points of $ \partial_{ B_r } \Sigma $, i.e., 
$ \hat{\sigma}^{ \ell ( w ) } \left( w, B_r, B_r, \ldots \right) = \left( B_r, B_r, \ldots \right) $.
Consequently, $ \hat{\mathcal{L}}_{0} $ absorbs $ \partial \Sigma $.

Note that the shift map $ \hat{\sigma} $ is measurable, since 
\begin{align*}  
\hat{\sigma}^{-1} \left( Z [ \epsilon ] \right)
& = \hat{\Sigma},
\\
\hat{\sigma}^{-1} \left( Z [ w ] \right) 
& = \bigsqcup_{ i \in \mathscr{P} (w) } Z [ i w ],
\\
\hat{\sigma}^{-1} \left( Z [ B_r ; S ] \right) 
& = \big\{ \left( B_r, B_r, \ldots \right) \big\} 
\, \sqcup \, 
\bigsqcup_{ i \in \mathscr{L}_1 } \,
Z [ i B_r; S ],
\\
\text{and} \qquad
\hat{\sigma}^{-1} \left( Z [ w B_r ; S ] \right) 
& = \bigsqcup_{ i \in \mathscr{P} (w) } \,
Z [ i w B_r; S ],
\end{align*}
where $ w $ is any finite non-empty word and $ 1 \leq r \leq s $.
The shift map $ \hat{\sigma} $ is continuous only on $ \hat{\Sigma} \setminus \hat{\mathcal{L}}_{0} $
(see Proposition~4.1 of \cite{AS:BSM21}).
In fact, the discontinuity of $ \hat{\sigma} $ on $ \hat{\mathcal{L}}_{0} $ can be observed from the convergences of the following sequences,
as $ k \in B_r $ tends to $ \infty $,
\begin{equation*}
\begin{array}{rl}
\left( k, 0, 0, \ldots \right)
& \to \,
\left( B_r, B_r, \ldots \right) 
\qquad \text{and} \\
\hat{\sigma} \left( k, 0, 0, \ldots \right) 
& \to \,
\left( 0, 0, \ldots \right) 
\neq \hat{\sigma} \left( B_r, B_r, \ldots \right) .
\end{array}
\end{equation*} 

\hypertarget{sec-invariant-probabilities}{}
\section{Invariant Probabilities}
\label{section-invariant-probabilities}

Let us introduce the following sets of invariant measures
\begin{align*}
\mathcal{M} \left( \Sigma, \sigma \right) 
& = \left\{ \mu \in \mathcal{M} \left( \Sigma \right) \,:\, 
\begin{array}{c}
\mu \text{ is $ \sigma $-invariant probability, i.e.,} \\
\sigma_{\ast} \mu ( \, \cdot \, ) = \mu \big( \sigma^{-1} ( \, \cdot \, ) \big) = \mu ( \, \cdot \, )
\end{array}    
\right\},
\\
\text{and} \quad
\mathcal{M} \big( \hat{\Sigma}, \hat{\sigma} \big) 
& = \left\{ \hat{\mu} \in \mathcal{M} \big( \hat{\Sigma} \big) \,:\,
\begin{array}{c}
\hat{\mu} \text{ is $ \hat{\sigma} $-invariant probability, i.e.,} \\
\hat{\sigma}_{\ast} \hat{\mu} ( \, \cdot \, ) = \hat{\mu} \big( \hat{\sigma}^{-1} ( \, \cdot \, ) \big) = \hat{\mu} ( \, \cdot \, ) 
\end{array} 
\right\}.
\end{align*}
Whenever we write  
$ \overline{ \mathcal{M} \left( \Sigma, \sigma \right) } $ or 
$ \overline{ \mathcal{M} \big( \hat{\Sigma}, \hat{\sigma} \big) } $, 
it will be understood that we consider the closure with respect to the weak$^{\ast}$ topology of $ \mathcal{M} \big( \hat{\Sigma} \big) $.

In this Section, we will study how these subsets are related, as well as their main properties.
It is clear that $ \mathcal{M} \left( \Sigma, \sigma \right) \subset \mathcal{M} \big( \hat{\Sigma}, \hat{\sigma} \big) $, and both sets are convex.
We cannot expect that neither $ \mathcal{M} \left( \Sigma, \sigma \right) $ nor $ \mathcal{M} \big( \hat{\Sigma}, \hat{\sigma} \big) $ to be closed in general, as ilustrated below.
\begin{example}  \label{example-invariant-mesures-not-closed}
From item~\ref{item-convegence-partial-sigma} of Lemma~\ref{lemma-convegence-hat-sigma}, we observe that the sequence of invariant probabilities
\begin{equation*}
\left\{ \frac{1}{2} \,\delta_{ \left( k, 0, k, 0, \ldots \right) } + \frac{1}{2} \, \delta_{ \left( 0, k, 0, k, \ldots \right) } \right\}_{ k \in B_r }
\subset \mathcal{M} \left( \Sigma, \sigma \right),
\end{equation*} 
as $ k \in B_r $ tends to $ \infty $, converges to
\begin{equation*}
\dfrac{1}{2} \, \delta_{ \left( B_r, B_r, \ldots \right) } + \dfrac{1}{2} \, \delta_{ \left( 0,  B_r, B_r, \ldots \right) }
\in \overline{ \mathcal{M} \left( \Sigma, \sigma \right) }
\subset \overline{ \mathcal{M} \big( \hat{\Sigma}, \hat{\sigma} \big) }.
\end{equation*}
However, this probability is not invariant:
\begin{equation*}
\hat{\sigma}_{\ast} \left( \dfrac{1}{2} \, \delta_{ \left( B_r, B_r, \ldots \right) } + \dfrac{1}{2} \, \delta_{ \left( 0,  B_r, B_r, \ldots \right) } \right) 
= \delta_{ \left( B_r, B_r, \ldots \right) }.
\end{equation*}
\end{example}

It is important to highlight that under suitable hypotheses on the original shift, precisely the \hyperlink{fcpa}{finite cyclic predecessor assumption} and 
the \hyperlink{dpm}{denseness of periodic measures}, we can guarantee the compactness of the set of blur invariant measures (see Corollary~\ref{corollary-blur-measure-compact}).
A first step will be a convex decomposition of blur invariant probabilities
(Proposition~\ref{proposition-decomposition-blur-invariant-measures}),
which will hold for any shift space over a countable alphabet.

\subsection{Convexity Description}
\label{subsection-convexity}

Our goal is to exhibit a convex characterization for 
$ \mathcal{M} \big( \hat{\Sigma}, \hat{\sigma} \big) $.
In order to do that, we need to better understand the behavior of blur shift probabilities on $ \partial \Sigma $.

\begin{lemma}  \label{lemma-behavior-blur-measure-partial-sigma} 
If $ \hat{\mu} \in \mathcal{M} \big( \hat{\Sigma}, \hat{\sigma} \big) $,
then 
$ \hat{\mu} \big( \, \bigsqcup_{k \geq 1} \, \hat{\mathcal{L}}_{k} \, \big) = 0 $.
In particular,
\begin{equation*}
\hat{\mu} \big( \hat{\mathcal{L}}_{0} \big) = \hat{\mu} \left( \partial \Sigma \right)
\qquad \text{and} \qquad
\hat{\mu} \left( \big\{ \left( B_r, B_r, \ldots \right) \big\} \right) = \hat{\mu} \left( \partial_{B_r} \Sigma \right).
\end{equation*}
\end{lemma}

\begin{proof}
\let\qed\relax
Remember that $ \hat{\mathcal{L}}_{0} $ absorbs $ \partial \Sigma $.
Since $ \hat{\mu} $ is $ \hat{\sigma} $-invariant measure, we have
$ \hat{\mu} \big( \hat{\mathcal{L}}_{0} \big) 
= \hat{\mu} \big( \hat{\sigma}^{-k} \big( \hat{\mathcal{L}}_{0} \big) \big)
= \hat{\mu} \big( \hat{\mathcal{L}}_{0} \sqcup \cdots \sqcup \hat{\mathcal{L}}_{k} \big) $
for all $ k \geq 1 $.
Hence,
\begin{equation*}
\hat{\mu} \bigg( \, \bigsqcup_{k \geq 1} \, \hat{\mathcal{L}}_{k} \, \bigg) 
= \lim_{ k \to \infty } \, \hat{\mu} \big( \hat{\mathcal{L}}_{0} \sqcup \cdots \sqcup \hat{\mathcal{L}}_{k} \big) - \hat{\mu} \big( \hat{\mathcal{L}}_{0} \big)
= 0.
\tag*{\qedsymbol}
\end{equation*}
\end{proof}

\begin{lemma}  \label{lemma-total-mass-partial-sigma}
If $ \hat{\nu} \in \mathcal{M} \big( \hat{\Sigma}, \hat{\sigma} \big) $ and $ \hat{\nu} ( \partial \Sigma ) = 1 $, 
then $ \hat{\nu} \in \mathrm{Conv} \left( \, \left\{ \delta_{ \left( B_r, B_r, \ldots \right) } \,:\, 1 \leq r \leq s  \right\} \, \right) $.
\end{lemma}

\begin{proof}
It is immediate from Lemma~\ref{lemma-behavior-blur-measure-partial-sigma} that
$ \hat{\nu} \big( \hat{\mathcal{L}}_{0} \big) = \hat{\nu} \left( \partial \Sigma \right) = 1 $.
Hence,
\begin{equation*}
\hat{\nu} \left( \, \cdot \, \right)
= \hat{\nu} \big( \, \cdot \, \cap \, \hat{\mathcal{L}}_{0} \big)
= \sum_{r=1}^{s} \, \hat{\nu} \left( \, \cdot \, \cap \big\{ \left( B_r, B_r, \ldots \right) \big\} \right)
= \sum_{r=1}^{s} \, \alpha_r \, \delta_{ \left( B_r, B_r, \ldots \right) }  \left( \, \cdot \, \right), 
\end{equation*}
where
$ \alpha_r : = \hat{\nu} \left( \big\{ \left( B_r, B_r, \ldots \right) \big\} \right) $ fulfills
$ \sum_{r=1}^{s} \alpha_r  = \hat{\nu} \big( \hat{\mathcal{L}}_{0} \big) = 1 $.
\end{proof}

\begin{figure}[ht]
\centering
\begin{tikzpicture}
\tikzstyle{linha} = [draw=black, line width=0.8pt]

\path[linha] (0,0)
to[out=-110,in=-80,looseness=0.9] (4,0);
\path[linha, draw=black!30] (4,0)
to[out=110,in=40,looseness=0.9] (0,0);
\node[] (Base) at (2,-0.15) {$ \mathcal{M} \left( \Sigma, \sigma \right) $};

\path[linha, black!70, dotted] (0.6,-0.2) -- (2,4);
\path[linha] (0,0) -- (2,4);
\path[linha] (4,0) -- (2,4);

\fill[black] (2,4) circle[radius=2pt];
\node[above] (Delta) at (2,4) {$ \delta_{ \left( B_1, B_1, \ldots \right) } $};

\fill[black] (0.6,-0.2) circle[radius=2pt];
\node[below] (Mu) at (0.6,-0.2) {$ \mu $};

\fill[black] (1.17,1.5) circle[radius=2pt];
\node[right] (Mu) at (1.17,1.5) {$ \hat{\mu} $};

\end{tikzpicture}
\qquad
\begin{tikzpicture}
\tikzstyle{linha} = [draw, line width=0.8pt]

\path[linha] (0,0)
to[out=-110,in=-80,looseness=0.9] (4,0);
\path[linha, draw=black!30] (4,0)
to[out=110,in=40,looseness=0.9] (0,0);
\node[] (Base) at (2,-0.15) {$ \mathcal{M} \left( \Sigma, \sigma \right) $};

\path[linha, draw=black!30] (0,0) -- (2.3,4);
\path[linha, draw=black!30] (4,0) -- (2.3,4);
\path[linha] (0,0) -- (1,3);
\path[linha] (4,0) -- (3,3.2);
\path[linha, black!70, dotted] (0.6,-0.2) -- (2.2,3.3);
\path[linha] (1,3) -- (3,3.2) -- (2.3,4) -- (1,3);

\fill[black] (1,3) circle[radius=2pt];
\node[left] (Delta1) at (1,3) {$ \delta_{ \left( B_1, B_1, \ldots \right) } $};
\fill[black] (2.3,4)  circle[radius=2pt];
\node[above] (Delta2) at (2.3,4) {$ \delta_{ \left( B_2, B_2, \ldots \right) } $};
\fill[black] (3,3.2) circle[radius=2pt];
\node[right] (Delta3) at (3,3.2) {$ \delta_{ \left( B_3, B_3, \ldots \right) } $};

\fill[black] (2.2,3.3) circle[radius=2pt];
\node[above] (Nu) at (2.2,3.3) {$ \hat{\nu} $};

\fill[black] (0.6,-0.2) circle[radius=2pt];
\node[below] (Mu) at (0.6,-0.2) {$ \mu $};

\fill[black] (1.38,1.5) circle[radius=2pt];
\node[right] (Mu) at (1.38,1.5) {$ \hat{\mu} $};

\end{tikzpicture}
\caption[Graphical representations of the set of $ \hat{\sigma} $-invariant measures]
{Graphical representations of convex decompositions of the set of $ \hat{\sigma} $-invariant measures for blur shifts with one-symbol and three-symbol resolutions}
\label{figure-decomposition-hat-invariant-measures}
\end{figure}
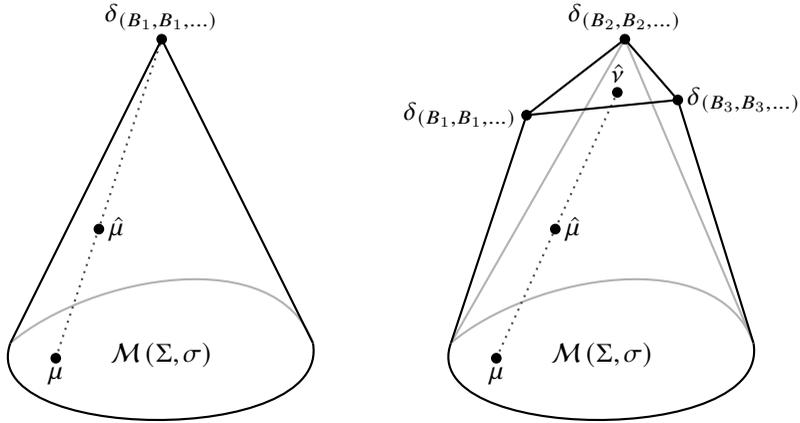

We obtain the following convex characterization for the $ \hat{\sigma} $-invariant probabilities
(see Figure~\ref{figure-decomposition-hat-invariant-measures}).
\begin{proposition}  \label{proposition-decomposition-blur-invariant-measures}
\enspace
$ \displaystyle
\mathcal{M} \big( \hat{\Sigma}, \hat{\sigma} \big) 
= \mathrm{Conv} \left( \, \mathcal{M} \left( \Sigma, \sigma \right) \, \sqcup \, \left\{ \delta_{ \left( B_r, B_r, \ldots \right) } \,:\, 1 \leq r \leq s  \right\} \, \right) $.
\end{proposition}

\begin{proof}
By the convexity of the set of blur invariant probabilities, it follows
$ \mathrm{Conv} \left( \, \mathcal{M} \left( \Sigma, \sigma \right) \, \cup \, \left\{ \delta_{ \left( B_r, B_r, \ldots \right) } \,:\, 1 \leq r \leq s \right\} \, \right)  \subset \mathcal{M} \big( \hat{\Sigma}, \hat{\sigma} \big) $.
On the other hand, assume that $ \hat{\mu} \in \mathcal{M} \big( \hat{\Sigma}, \hat{\sigma} \big) $ and remember that $ \hat{\Sigma} = \Sigma \sqcup \partial \Sigma $.
We will consider three cases:

\smallskip\noindent
{\itshape Case $ \hat{\mu} \left( \Sigma \right) = 0 $.}
It is immediate that $ \hat{\mu} \left( \partial \Sigma \right) = 1 $.
Thus, Lemma~\ref{lemma-total-mass-partial-sigma} guarantees that $ \hat{\mu} \in \mathrm{Conv} \left( \, \left\{ \delta_{ \left( B_r, B_r, \ldots \right) } \,:\, 1 \leq r \leq s  \right\} \, \right) $.

\smallskip\noindent
{\itshape Case $ \hat{\mu} (\Sigma) \in (0, 1) $.} 
Consider the decomposition
\begin{equation*}
\hat{\mu} ( \,\cdot\, ) 
= \hat{\mu} ( \, \cdot \, \cap \, \Sigma ) + \hat{\mu} ( \, \cdot \, \cap \, \partial \Sigma ) 
= t \, \dot{\mu} ( \,\cdot\, ) + ( 1 - t ) \, \hat{\nu} ( \, \cdot \, ),
\end{equation*}
with constant
$ t = \hat{\mu} (\Sigma) $,
consequently $ 1 - t = \hat{\mu} \big( \partial \Sigma \big) $,
and $ \hat{\sigma} $-invariant probabilities
\begin{equation*}
\dot{\mu} ( \,\cdot\, ) = \dfrac{ \hat{\mu} \big( \cdot \, \cap \, \Sigma \big) }{ \hat{\mu} (\Sigma) }
\qquad \text{and} \qquad
\hat{\nu} ( \,\cdot\, ) = \dfrac{ \hat{\mu} \big( \cdot \, \cap \, \partial \Sigma \big) }{ \hat{\mu} \big( \partial \Sigma \big) }.
\end{equation*}
In particular, $ \dot{\mu} \left( \Sigma \right) = 1 $ and $ \hat{\nu} \left( \partial \Sigma \right) = 1 $.
By Lemma~\ref{lemma-topological-inclusion}, we can identify $ \dot{\mu} $ as a probability measure in $ \mathcal{M} \left( \Sigma, \sigma \right) $.
Moreover, 
$ \nu = \sum_{r=1}^{s} \alpha_r \, \delta_{ \left( B_r, B_r, \ldots \right) } $ with $ \sum_{r=1}^{s} \alpha_r = 1 $, due to Lemma~\ref{lemma-total-mass-partial-sigma}.
We conclude that
\begin{equation*}
\hat{\mu} 
= t \, \dot{\mu} + \sum_{r=1}^{n}  ( 1 - t ) \, \alpha_r \, \delta_{ \left( B_r, B_r, \ldots \right) } 
\quad\text{with}\quad
t + \sum_{r=1}^{n} ( 1 - t ) \, \alpha_r = 1.
\end{equation*}

\smallskip\noindent
{\itshape Case $ \hat{\mu} (\Sigma) = 1 $.}
It is straightforward from Lemma~\ref{lemma-topological-inclusion}.
\end{proof}

\begin{remark}
From the previous result, for any $ \hat{\mu} \in \mathcal{M} \big( \hat{\Sigma}, \hat{\sigma} \big) $,
we have 
$ \hat{\mu} \left( \, \cdot \, \cap \Sigma \right) 
= t \, \dot{\mu} \left( \, \cdot \, \right) $
with $ t = \hat{\mu} \left( \Sigma \right) \in [0,1] $ 
and $ \dot{\mu} \in \mathcal{M} \left( \Sigma, \sigma \right) $.
This fact induces a projection from the set of $ \hat{\sigma} $-invariant probabilities on $ \hat{\Sigma} $ to the set of $ \sigma $-invariant subprobabilities on $ \Sigma $, which is a central notion used in~\cite{IV:JAM21} to study invariant measures for countable Markov shifts.
The phenomenon of escape of mass can be translated as a transference of mass from $ \Sigma $ to points of $ \partial \Sigma $. 
\end{remark}

\subsection{Compactness}
\label{subsection-compactness}
In this subsection, we show the compactness of the set of blur invariant probabilities, supposing 
the \hyperlink{fcpa}{finite cyclic predecessor assumption} and the
\hyperlink{dpm}{denseness of periodic measures}.
Initially, we argue that these hypotheses on the original shift impose that any measure belonging to the closure of $ \sigma $-invariant probabilities fulfills the necessary condition to be a blur invariant probability described in Lemma~\ref{lemma-behavior-blur-measure-partial-sigma}, i.e., the set $ \partial \Sigma \setminus \hat{\mathcal{L}}_0 $ has null measure.

\begin{lemma} \label{lemma-fcpa-dpm-l1-null-measure}
Let $ \Sigma $ be a shift over a countable alphabet verifying 
the \hyperlink{fcpa}{finite cyclic predecessor assumption} and 
the \hyperlink{dpm}{denseness of periodic measures}.
Then,
\begin{equation*}
\bar{\mu} \big( \hat{\mathcal{L}}_1 \big) = 0 
\quad \text{for every} \quad
\bar{\mu} \in \overline{ \mathcal{M} \left( \Sigma, \sigma \right) }.
\end{equation*}
\end{lemma}

\begin{proof}
\let\qed\relax
By Lemma~\ref{lemma-topological-inclusion} and by the \hyperlink{dpm}{denseness of periodic measures}, the set of ergodic probabilities supported on periodic orbits of $ \Sigma $ is also (weak$^{\ast}$) dense on $ \overline{ \mathcal{M} \left( \Sigma, \sigma \right) } $ as a subset of $ \mathcal{M} \big( \hat{\Sigma} \big) $.
Therefore, given $ \bar{\mu} \in \overline{ \mathcal{M} \left( \Sigma, \sigma \right) } $, let 
\begin{equation*}
y^l = \big( y_{0}^{l}, \ldots, y_{p_{l}-1}^{l}, y_{0}^{l}, \ldots \big),
\quad \text{where} \enspace
p_{l} = \min \left\{ p \geq 1 \,:\, 
y^{l} = \sigma^{p} \big( y^{l} \big) \right\},
\quad \text{for} \enspace l \geq 1,
\end{equation*}
denote a sequence of periodic points such that the respective periodic probabilities 
$ \dot{\mu}_l = \frac{1}{ p_{l} } \sum_{ k = 0 }^{ p_{l} - 1 } \delta_{ \sigma ( y^{l} ) } $
converge to $ \bar{\mu} $.

Note that the countable set 
$ \hat{\mathcal{L}}_1 
= \bigsqcup_{ 1 \leq r \leq s } \,
\bigsqcup_{ j \in B_r }
\bigsqcup_{ i \in \mathscr{P} (j) }
\big\{ \left( i, B_r, B_r, \ldots \right) \big\} $
has null measure if, and only if, 
$ \bar{\mu} \left( \big\{ \left( a, B_r, B_r, \ldots \right) \big\} \right) = 0 $
for all $ a \in \bigsqcup_{ j \in B_r } \mathscr{P} (j) $ and for any $ 1 \leq r \leq s $.
We fix once for all $ \left( a, B_r, B_r, \ldots \right) $ and we consider the following three possibilities:

\smallskip \noindent
{\itshape --
Suppose that $ \big\{ \max_{ 0 \leq j \leq p_l} \, y_{j}^{l} \big\}_{ l \geq 1 } $ is bounded}. 
Take into account the finite subset of $ B_r $ given as
\begin{equation*}
S = B_r \cap \left\{ 1, \ldots, \, \max_{ l \geq 1 } \, \max_{ 0 \leq j \leq p_l} \, y_{j}^{l} \, \right\},
\end{equation*}
and consider the clopen neighborhood $ Z [ a B_r ; S ] $ of $ \left( a, B_r, B_r, \ldots \right) $.
Since $ y_{j}^{l} \notin B_r \setminus S $, for every $ 0 \leq j \leq p_l $,
it follows that 
$ \mathrm{supp} \, \dot{\mu}_l 
= \mathrm{orb} \left( y^l \right)
\subset Z [ a B_r ; S ]^{\complement} $.
In particular, $ \dot{\mu}_{l} \left( Z [ a B_r; S ] \right) = 0 $ for all $ l \geq 1 $.
Due to Portmanteau Theorem, we have that
\begin{equation*}
\bar{\mu} \left( \big\{ \left( a, B_r, B_r, \ldots \right) \big\} \right)
\leq \bar{\mu} \left( Z [ a B_r; S ] \right) 
= \lim_{ l \to \infty } \, \dot{\mu}_{l} \left( Z [ a B_r; S ] \right)
= 0.
\end{equation*}

\smallskip \noindent
{\itshape --
Suppose that 
$ \big\{ \max_{ 0 \leq j \leq p_l} y_{j}^{l} \big\}_{ l \geq 1 } $ is unbounded 
and $ \left\{ p_{l} \right\}_{ l \geq 1 } $ is bounded}. 
By taking a subsequence, we may assume that $ p_{l'} = p $ is constant.
Without loss of generality, we consider $ \big\{ y_{0}^{l'} \big\} $ unbounded.
Note that
\begin{equation*}
\big\{ y_{0}^{l'} \big\} 
\subset
\bigcup_{ l' \geq 1 } \, \mathscr{P} ( y_{1}^{l'} ) \cap \mathscr{F}_{ p - 1 } ( y_{1}^{l'} ).
\end{equation*}
By the \hyperlink{fcpa}{finite cyclic predecessor assumption}, we have a union of finite sets that contains an infinite set.
This is only possible if $ \big\{ y_{1}^{l'} \big\} $ is also a unbounded set.
By a recursive argument, we obtain that $ \big\{ y_{j}^{l'} \big\} $ is unbounded for any $ 0 \leq j \leq p-1 $.
In particular, there exists an infinite subset of indexes $ l'' $ such that $ y_{j}^{l''} > a $ for all $ 0 \leq j \leq p - 1 $.
We thus obtain
$ \mathrm{supp} \, \dot{\mu}_{l''} 
= \mathrm{orb} \left( y^{l''} \right)
\subset Z [ a B_r ]^{\complement} $,
so that
\begin{equation*}
\bar{\mu} \left( \big\{ \left( a, B_r, B_r, \ldots \right) \big\} \right)
\leq \bar{\mu} \left( Z [ a B_r ] \right) 
= \lim_{ l'' \to \infty } \, \dot{\mu}_{l''} \left( Z [ a B_r ] \right) = 0.
\end{equation*}

\smallskip \noindent
{\itshape --
Suppose that 
$ \big\{ \max_{ 0 \leq j \leq p_l} y_{j}^{l} \big\}_{ l \geq 1 } $ 
and $ \left\{ p_{l} \right\}_{ l \geq 1 } $ are unbounded}. 
Given a positive integer $ m $, consider the following finite sets
\begin{align*}
R_{0} & = \left\{a \right\}, 
& 
R_{1} & = \bigcup_{ k = 1 }^{m} \,
\mathscr{P} (a) \cap \mathscr{F}_{k} (a),
\\
R_{2} &= \bigcup_{ b \in R_{1} } \, \bigcup_{ k = 1 }^{m} \,
\mathscr{P} (b) \cap \mathscr{F}_{k} (b),
& \ldots, \qquad
R_{m} &= \bigcup_{ c \in R_{m-1} } \, \bigcup_{ k = 1 }^{m} \,
\mathscr{P} (c) \cap \mathscr{F}_{k} (c).
\end{align*}
Denote $ S_m = B_r \cap \left( R_{0} \cup R_{1} \cup R_{2} \cup \cdots \cup R_{m} \right) $.
Due to Portmanteau Theorem, we have that
\begin{equation*}
\bar{\mu} \left( \big\{ \left( a, B_r, B_r, \ldots \right) \big\} \right)
\leq \bar{\mu} \left( Z [ a B_r; S_m ] \right) 
= \lim_{ l \to \infty } \, \dot{\mu}_{l} \left( Z [ a B_r; S_m ] \right).
\end{equation*}
To conclude, we will show that
$ \dot{\mu}_{l''} \left( Z [ a B_r; S_m ] \right) < \frac{1}{ m + 1 } $
for infinitely many indexes~$ l'' $.
Without lost of generality, let $ \left\{ y_{0}^{l} \right\} $ be unbounded.
We may consider a subsequence such that
$ p_{l'} > m + 1 $, for all $ l' $.
From the \hyperlink{fcpa}{finite cyclic predecessor assumption} and a recursive argument as above, we obtain that
$ \big\{ y_{j}^{l'} \big\} $ is also unbounded for any $ 0 \leq j \leq m $.
Passing to another subsequence, we may assume that $ y_{j}^{l''} > a $, for all $ l'' $ and $ 0 \leq j \leq m $.
Note that:
\begin{itemize}
\item $ \sigma^i \left( y^{l''} \right) \in Z [ a B_r ; S_m ] $ 
implies $ m < i < p_{l} $, 
because $ y_{j}^{l''} > a $ for all $ 0 \leq j \leq m $;
\item $ \sigma^i \left( y^{l''} \right) $, $ \sigma^j \left( y^{l''} \right) \in Z [ a B_r ; S_m ] $, with $ i < j $,
implies $ j - i > m + 1 $, 
because, otherwise, if $ 1 < j - i = k \leq m + 1 $, as
$ y_{i}^{l''} = y_{j}^{l''} = a $ 
and $ y_{i+1}^{l''} \notin R_{k-1} $,
then
$ y_{j-1}^{l''} \in \mathscr{P} (a) $ and $ y_{j-1}^{l''} \notin R_{1} $ 
which means that 
$ y_{j-1}^{l''} \notin \mathscr{F}_{k-1} (a) $, which is a contradiction.
\end{itemize}
Hence, from a counting argument, it follows that
\begin{equation*}
\dot{\mu}_{l''} \left( Z [ a B_r ; S_n ] \right)  
< \dfrac{1}{ p_{l''} } \left\lceil \dfrac{ p_{l''} - m - 1 }{ m + 1 } \right\rceil
\leq \dfrac{1}{ p_{l''} } \left( \dfrac{ p_{l''} - m - 1 }{ m + 1 } + 1 \right)
= \dfrac{ 1 }{ m + 1 },
\quad \text{for all} \enspace l''.
\tag*{\qedsymbol}
\end{equation*}
\end{proof}

\begin{remark}
For countable Markov shifts, we point out the condition $ \bar{\mu} \big( \hat{\mathcal{L}}_1 \big) = 0 $ for all $ \bar{\mu} \in \overline{ \mathcal{M} \left( \Sigma, \sigma \right) } $ implies the \hyperlink{fcpa}{finite cyclic predecessor assumption}.
In fact, assume that $ \mathscr{P} (a) \cap \mathscr{F}_{m} (a) $ is infinite, for some $ a \in \mathscr{L}_1 $ and $ m \geq 1 $.
For each $ l \in \mathscr{P} (a) \cap \mathscr{F}_{m} (a) $, thus there exists a periodic point $ y^l = \left( a, u^{l}, l, a, u^{l}, l \ldots \right) \in \Sigma $, where $ u^{l} \in \mathscr{L}_{m-1} $ is a word that connects $ a $ to $ l $.
Because $ \mathcal{M} \big( \hat{\Sigma} \big) $ is compact, there is a subsequence of periodic measures
$ \dot{\mu}_{l'} = \frac{1}{ m + 1 } \, \sum_{ i = 0 }^{m} \delta_{ \sigma^i ( y^{l'} ) } $ 
converging to some $ \bar{\mu} \in \overline{ \mathcal{M} \left( \Sigma, \sigma \right) } $, as $ l' \in \mathscr{P} (a) \cap \mathscr{F}_{m} (a) $ tends to $ \infty $.
Note that there is a subsequence $ y^{l''} $ converging to  a point $ \left( a, v, B_r, B_r, \ldots \right) \in \hat{\Sigma} $, with $ \ell(v) < m $ and $ 1 \leq r \leq s $.
For any clopen neighborhood $ Z [ a v B_r ; S ] $ of $ \left( a, v, B_r, B_r, \ldots \right) $, it is easy to see that
$ \bar{\mu}_{l''} \left( Z [ a v B_r; S ] \right) 
\geq \frac{1}{ m + 1 } \, \delta_{ y^{l''} } \left( Z [ a v B_r; S ] \right) 
= \frac{1}{ m + 1 } $ for $ l'' $ sufficiently large.
By the Portmanteau Theorem, 
considering $ S_n = B_r \cap \left\{ 1, \ldots, n \right\} $,
we obtain that
\begin{align*}
\bar{\mu} \big( \hat{\mathcal{L}}_{\ell(v)+1} \big) 
\geq \bar{\mu} \left( \big\{ \left( a, v, B_r, B_r, \ldots \right) \big\} \right)
= \inf_{n \geq 1} \, \lim_{ l'' \to \infty } \, \dot{\mu}_{l''} \left( Z [ a v B_r; S_n ] \right)
\geq \dfrac{1}{ m + 1 } 
> 0.
\end{align*}
For $ \ell(v) \geq 1 $, the contrapositive statement is a consequence of the following lemma for general shifts. 
\end{remark}

\begin{lemma}  \label{lemma-lk-decrease-measure}
Let $ \Sigma $ be a shift over a countable alphabet.
If $ \bar{\mu} \in \overline{ \mathcal{M} \left( \Sigma, \sigma \right) } $,
then $ \bar{\mu} \big( \hat{\mathcal{L}}_{k+1} \big) \leq \bar{\mu} \big( \hat{\mathcal{L}}_{k} \big) $ for all $ k \geq 1 $.
\end{lemma}

\begin{proof}
Let $ \{ \mu_n \} \in \mathcal{M} \big( \Sigma, \sigma \big) 
\subset \mathcal{M} \big( \hat{\Sigma}, \hat{\sigma} \big) $ be a sequence that converges to $ \bar{\mu} \in \mathcal{M} \big( \hat{\Sigma} \big) $.
Fix a point $ \left( w, B_r, B_r, \ldots \right) $, with $ \ell (w) \geq 1 $, and consider any finite subset $ S $ of $ B_r $. 
The clopen $ Z [ w B_r ; S ] $ that does not intersect $ \hat{\mathcal{L}}_0 $.
Since $ \hat{\sigma} $ is continuous outside $ \hat{\mathcal{L}}_0 $, the inverse image
$ \hat{\sigma}^{-1} \left( Z [ w ; B_r ; S ] \right) $
is an open neighborhood of 
$ \hat{\sigma}^{-1} \left( \big\{ \left( w, B_r, B_r, \ldots \right) \big\} \right) $.
Due to Portmanteau Theorem and $ \hat{\sigma} $-invariance of $ \mu_n $, we obtain 
\begin{align*}
\bar{\mu} \big( \hat{\sigma}^{-1} \left( \big\{ \left( w, B_r, B_r, \ldots \right) \big\} \right) \big)
& \leq \bar{\mu} \big( \hat{\sigma}^{-1} \left( Z [ w B_r ; S  ] \right) \big) 
\leq \liminf_{ n \to \infty } \, \mu_n  \big( \hat{\sigma}^{-1} \left( Z [ w B_r ; S  ] \right) \big) \\
& = \lim_{ n \to \infty } \, \mu_n \left( Z [ w B_r ; S  ] \right) 
= \bar{\mu} \left( Z [ w B_r ; S  ] \right).
\end{align*}
As this inequality holds for every finite subset $ S $ of $ B_r $, we have that
$ \bar{\mu} \big( \hat{\sigma}^{-1} \left( \big\{ \left( w, B_r, B_r, \ldots \right) \big\} \right) \big) 
\leq \bar{\mu} \left( \big\{ \left( w, B_r, B_r, \ldots \right) \big\} \right) $.
To conclude the proof, note that
$ \hat{\mathcal{L}}_{k} $ is a countable union of singletons $ \big\{ \left( w, B_r, B_r, \ldots \right) \big\} $, with $ \ell (w) = k $ and $1 \leq r \leq s$,
and consequently $ \hat{\mathcal{L}}_{k+1} = \hat{\sigma}^{-1} \big( \hat{\mathcal{L}}_{k} \big) $ is countable disjoint union of inverse images 
$ \hat{\sigma}^{-1} \left( \big\{ \left( w, B_r, B_r, \ldots \right) \big\} \right) $.
\end{proof}

In the next result, we provide equivalent formulations for the compactness of the set of blur invariant probability measures.

\begin{proposition}  \label{proposition-equivalent-blur-measure-compact}
Let $ \Sigma $ be a shift over a countable alphabet.
The following statements are equivalent:
\begin{enumerate}[label=(\roman*)]
\setlength\itemsep{4pt}
\item\label{item-l1-null-measure}
$ \bar{\mu} \big( \hat{\mathcal{L}}_1 \big) = 0 $, 
for every $ \bar{\mu} \in \overline{ \mathcal{M} \left( \Sigma, \sigma \right) } $;
\item\label{item-sum-lk-null-measure}
$ \bar{\mu} \left( \bigsqcup_{ k \geq 1 } \hat{\mathcal{L}}_k \right) = 0 $,
for all $ \bar{\mu} \in \overline{ \mathcal{M} \left( \Sigma, \sigma \right) } $;
\item\label{item-closure-subset-blur-invariant-measure} 
$ \overline{ \mathcal{M} \left( \Sigma, \sigma \right) } \subset \mathcal{M} \big( \hat{\Sigma}, \hat{\sigma} \big) $;
\item\label{item-blur-invariant-measure-compact} 
$ \mathcal{M} \big( \hat{\Sigma}, \hat{\sigma} \big) $ is (weak$^{\ast}$) compact.
\end{enumerate}
\end{proposition}

\begin{proof}
The equivalence
\ref{item-l1-null-measure}~$ \Leftrightarrow $~\ref{item-sum-lk-null-measure} 
is a straightforward consequence of Lemma~\ref{lemma-lk-decrease-measure}.
The implication 
\ref{item-closure-subset-blur-invariant-measure}~$ \Rightarrow $~\ref{item-sum-lk-null-measure} 
is immediate from Lemma~\ref{lemma-behavior-blur-measure-partial-sigma}.
Moreover,
\ref{item-blur-invariant-measure-compact}~$ \Rightarrow $~\ref{item-closure-subset-blur-invariant-measure}
because, when assuming \ref{item-blur-invariant-measure-compact},
$ \mathcal{M} \big( \hat{\Sigma}, \hat{\sigma} \big) $ is closed.
It remains to prove that
\ref{item-sum-lk-null-measure}~$ \Rightarrow $~\ref{item-closure-subset-blur-invariant-measure}~$ \Rightarrow $~\ref{item-blur-invariant-measure-compact}.

\smallskip
{\itshape \ref{item-sum-lk-null-measure}~$ \Rightarrow $~\ref{item-closure-subset-blur-invariant-measure}.}
Let $ \{ \mu_n \} \in \mathcal{M} \big( \Sigma, \sigma \big) \subset \mathcal{M} \big( \hat{\Sigma}, \hat{\sigma} \big) $ be a sequence that converges to a probability $ \bar{\mu} \in \overline{ \mathcal{M} \left( \Sigma, \sigma \right) } $.
In order to show that $ \bar{\mu} \in \mathcal{M} \big( \hat{\Sigma}, \hat{\sigma} \big) $, it is sufficient to verify the $ \hat{\sigma} $-invariance of the measure on generalized cylinder sets.

\smallskip\noindent
{\itshape Case $ Z [ \epsilon ] = \hat{\Sigma} $}. It is trivial.

\smallskip\noindent
{\itshape Case $ Z [ u ] $ with $ \ell (u) \geq 1 $}.
Remember that $ Z [ u ] $ is a the clopen set that does not intersect $ \hat{\mathcal{L}}_0 $.
Since $ \hat{\sigma} $ is continuous outside $ \hat{\mathcal{L}}_0 $, the inverse image
$ \hat{\sigma}^{-1} \left( Z [ u ] \right) $
is open.
Due to Portmanteau Theorem and $ \hat{\sigma} $-invariance of $ \mu_n $, we obtain 
\begin{equation}  \label{equation-semi-invariance}
\bar{\mu} \big( \hat{\sigma}^{-1} \left( Z [ u ] \right) \big) 
\leq \liminf_{ n \to \infty } \, \mu_n  \big( \hat{\sigma}^{-1} \left( Z [ u ] \right) \big) \\
= \lim_{ n \to \infty } \, \mu_n \left( Z [ u ] \right) 
= \bar{\mu} \left( Z [ u ] \right).
\end{equation}
Note that, for each $ k \geq 1 $, we have the partitions
\begin{align*}
\hat{\Sigma} \setminus \big( \hat{\mathcal{L}}_{0} \sqcup \cdots \sqcup \hat{\mathcal{L}}_{k-1} \big)
& = \bigsqcup_{ v \in \mathscr{L}_k } \, Z [ v ] 
\quad \text{and} \\
\hat{\Sigma} \setminus \big( \hat{\mathcal{L}}_{0}  \sqcup \cdots \sqcup \hat{\mathcal{L}}_{k} \big) 
= \hat{\Sigma} \setminus \hat{\sigma}^{-1} \big( \hat{\mathcal{L}}_{0} \sqcup \cdots \sqcup \hat{\mathcal{L}}_{k-1} \big)
& = \bigsqcup_{ v \in \mathscr{L}_k } \hat{\sigma}^{-1} \left( Z [ v ] \right).
\end{align*}
From item~\ref{item-sum-lk-null-measure}, it is easy to see that 
\begin{equation*}
\bar{\mu} \, \left( \hat{\Sigma} \setminus \big( \hat{\mathcal{L}}_{0} \sqcup \cdots \sqcup \hat{\mathcal{L}}_{k-1} \big) \right) 
= \bar{\mu} \, \big( \hat{\Sigma} \setminus \hat{\mathcal{L}}_0 \big) 
= \bar{\mu} \, \left( \hat{\Sigma} \setminus \big( \hat{\mathcal{L}}_{0} \sqcup \cdots \sqcup \hat{\mathcal{L}}_{k} \big) \right).
\end{equation*}
Hence, for any $ u $ with $ \ell(u) = k  \geq 1 $, the inequality~\eqref{equation-semi-invariance} allows us to obtain its opposite as follows:
\begin{align*}
\bar{\mu} \left( Z [ u ] \right)
& = \bar{\mu} \, \big( \hat{\Sigma} \setminus \hat{\mathcal{L}}_0 \big)
- \sum_{ v \in \mathscr{L}_k \setminus \{ u \} } \bar{\mu} \left( Z [ v ] \right) \\
& \leq \bar{\mu} \, \big( \hat{\Sigma} \setminus \hat{\mathcal{L}}_0 \big)
- \sum_{ v \in \mathscr{L}_k \setminus \{ u \} } \bar{\mu} \left( \hat{\sigma}^{-1} \left( Z [ v ] \right) \right) 
= \bar{\mu} \, \left( \hat{\sigma}^{-1} \left( Z [ u ] \right) \right).
\end{align*}

\smallskip\noindent
{\itshape Case $  Z [ B_r ; S ] $}. 
First, observe that
\begin{align*}
Z [ B_r ; S ] 
= \big\{ \left( B_r, B_r, \ldots \right) \big\}
& \, \sqcup \,
\bigsqcup_{ j \in B_r \setminus S } \, 
Z [ j ] 
\quad \text{and} \\
\hat{\sigma}^{-1} \left( Z [ B_r ; S ] \right) 
= \hat{\sigma}^{-1} \left( \big\{ \left( B_r, B_r, \ldots \right) \big\} \right)
&\, \sqcup \, 
\bigsqcup_{ j \in B_r \setminus S } \, \hat{\sigma}^{-1} \left( Z [ j ] \right),
\end{align*}
as well as
\begin{equation*}
\hat{\sigma}^{-1} \left( \big\{ \left( B_r, B_r, \ldots \right) \big\} \right)
= \big\{ \left( B_r, B_r, \ldots \right) \big\} 
\, \sqcup \,
\big\{ \left( i, B_r, B_r, \ldots \right) \,:\, \mathscr{F}_{1} (i) \cap B_r \text{ is infinite} \, \big\}    
\end{equation*}
Due to item~\ref{item-sum-lk-null-measure}, the measurable set
$ \left\{ \left( i, B_r, B_r, \ldots \right) \,:\, \mathscr{F}_{1} (i) \cap B_r \text{ is infinite} \right\} 
\subset \hat{\mathcal{L}}_1 $
has null measure and, therefore, 
$ \bar{\mu} \left( \hat{\sigma}^{-1} \left( \big\{ \left( B_r, B_r, \ldots \right) \big\} \right) \right) = \bar{\mu} \left(  \big\{ \left( B_r, B_r, \ldots \right) \big\} \right) $.
Thanks to the previous case, we have
$ \bar{\mu} \left( Z [ j ] \right) = \bar{\mu} \big( \hat{\sigma}^{-1} \left( Z [ j ] \right) \big) $ for every $ j \in B_r \setminus S $,
which shows that $ \bar{\mu} \left( Z [ B_r ; S ] \right) = \bar{\mu} \big( \hat{\sigma}^{-1} \left( Z [ B_r ; S ] \right) \big) $.

\smallskip\noindent
{\itshape Case $  Z [ w B_r ; S ] $ with $ \ell (w) \geq 1 $}. 
Note that
\begin{align*}
Z [ w B_r ; S ] 
& = \big\{ \left( w, B_r, B_r, \ldots \right) \big\} \, \sqcup \,  
\bigsqcup_{ j \in B_r \setminus S } Z [ w j ] 
\\
\hat{\sigma}^{-1} \left( Z [ w B_r ; S ] \right) 
& = \bigsqcup_{ i \in \mathscr{P} \left( w \right) } \big\{ \left( i, w, B_r, B_r, \ldots \right) \big\} \, \sqcup \, 
\bigsqcup_{ j \in B_r \setminus S } \,
\hat{\sigma}^{-1} \left( Z [ w j ] \right).
\end{align*}
From item~\ref{item-sum-lk-null-measure}, both
$ \big\{ \left( w, B_r, B_r, \ldots \right) \big\} \subset \hat{\mathcal{L}}_{ \ell(w) } $ and
$ \bigsqcup_{ i \in \mathscr{P} \left(  w \right) } \big\{ \left( i, w, B_r, B_r, \ldots \right) \big\} \subset \hat{\mathcal{L}}_{ \ell(w) + 1 } $ have null measure.
As we already know that
$ \bar{\mu} \left( Z [ w j ] \right) = \bar{\mu} \big( \hat{\sigma}^{-1} \left( Z [ w j ] \right) \big) $ for all $ j \in B_r \setminus S $,
we conclude that
$ \bar{\mu} \left( Z [ w B_r ; S ] \right) = \bar{\mu} \big( \hat{\sigma}^{-1} \left( Z [ w B_r ; S ] \right) \big) $.

\smallskip
{\itshape \ref{item-closure-subset-blur-invariant-measure}~$ \Rightarrow $~\ref{item-blur-invariant-measure-compact}.}
Consider a sequence $ \{ \hat{\mu}_n \} \in \mathcal{M} \big( \hat{\Sigma}, \hat{\sigma} \big) $.
By the convex decomposition of Proposition~\ref{proposition-decomposition-blur-invariant-measures}, we may write 
\begin{equation*}
\hat{\mu}_n = t_n \, \dot{\mu}_n + ( 1 - t_n ) \, \hat{\nu}_n,
\end{equation*}
where
$ t_n \in [ 0, 1 ] $,
$ \dot{\mu}_n \in \mathcal{M} \left( \Sigma, \sigma \right) $
and
$ \hat{\nu}_n \in \mathrm{Conv} \left( \, \left\{ \delta_{ \left( B_r, B_r, \ldots \right) } \,:\, 1 \leq r \leq s  \right\} \, \right) $.
Note that:
\begin{itemize}
\item there exists a subsequence $ \left\{ t_{n'} \right\} $ whose limit is $ t \in [0,1] $;
\item by compactness of the simplex $ \mathrm{Conv} \left( \, \left\{ \delta_{ \left( B_r, B_r, \ldots \right) } \,:\, 1 \leq r \leq s \right\} \, \right) $, there exists a subsequence $ \left\{ \hat{\nu}_{n''} \right\} $ that converges to $ \bar{\nu} \in \mathrm{Conv} \left( \, \left\{ \delta_{ \left( B_r, B_r, \ldots \right) } \,:\, 1 \leq r \leq s  \right\} \, \right) $;
\item since $ \overline{ \mathcal{M} \left( \Sigma, \sigma \right) } \subset \mathcal{M} \big( \hat{\Sigma} \big) $ is also compact, there exists a convergent subsequence $ \left\{ \dot{\mu}_{n'''} \right\} $ such that, 
due to the item~\ref{item-closure-subset-blur-invariant-measure} and the Proposition~\ref{proposition-decomposition-blur-invariant-measures}, its limit probability is given as 
$ \tau \, \dot{\mu} + ( 1 - \tau ) \, \dot{\nu} \in \mathcal{M} \big( \hat{\Sigma}, \hat{\sigma} \big) $, with
$ \tau \in [ 0, 1 ] $,
$ \dot{\mu} \in \mathcal{M} \left( \Sigma, \sigma \right) $ and
$ \dot{\nu} \in \mathrm{Conv} \left( \, \left\{ \delta_{ \left( B_r, B_r, \ldots \right) } \,:\, 1 \leq r \leq s  \right\} \, \right) $.
\end{itemize}
Taking this last subsequence, we obtain that
\begin{equation*}
\hat{\mu}_{n'''} 
\, \stackrel{\ast}{\rightharpoonup} \,
\bar{\mu} 
= t \, \big[ \tau \, \dot{\mu} + ( 1 - \tau ) \, \dot{\nu} \big]  + ( 1 - t ) \, \bar{\nu}.
\end{equation*}
If $ t = \tau = 1 $, we have
$ \bar{\mu} 
= \dot{\mu} \in \mathcal{M} \left( \Sigma, \sigma \right) 
\subset \mathcal{M} \big( \hat{\Sigma}, \hat{\sigma} \big) $.
Otherwise, we may write
\begin{equation*}
\bar{\mu} 
= ( t \tau ) \, \dot{\mu}  + ( 1 - t \tau ) \left[ \dfrac{ ( t - t \tau ) }{ ( 1 - t \tau ) } \, \dot{\nu} + \dfrac{ ( 1 - t ) }{ ( 1 - t \tau ) } \, \bar{\nu} \right],
\end{equation*}
with
$ t \tau \in [ 0, 1 ) $,
$ \dot{\mu} \in \mathcal{M} \left( \Sigma, \sigma \right) $ and
$ \frac{ ( t - t \tau ) }{ ( 1 - t \tau ) } \, \dot{\nu} + \frac{ ( 1 - t ) }{ ( 1 - t \tau ) } \, \bar{\nu} \in \mathrm{Conv} \left( \, \left\{ \delta_{ \left( B_r, B_r, \ldots \right) } \,:\, 1 \leq r \leq s  \right\} \, \right) $.
Hence, by Proposition~\ref{proposition-decomposition-blur-invariant-measures}, $ \bar{\mu} \in \mathcal{M} \big( \hat{\Sigma}, \hat{\sigma} \big) $.
This shows that $ \mathcal{M} \big( \hat{\Sigma}, \hat{\sigma} \big) $ is (sequentially) compact.
\end{proof}

We can close the discussion of this subsection by achieving as a consequence of the previous proposition, in view of Lemma~\ref{lemma-fcpa-dpm-l1-null-measure}, the compactness of the set of blur invariant probabilities.
\begin{corollary}  \label{corollary-blur-measure-compact}
Let $ \Sigma $ be a shift over a countable alphabet that satisfies both
\hyperlink{fcpa}{finite cyclic predecessor assumption}
and
\hyperlink{dpm}{denseness of periodic measures}.
Then, $ \mathcal{M} \big( \hat{\Sigma}, \hat{\sigma} \big) $ is (weak$^{\ast}$) compact.
\end{corollary}

\hypertarget{sec-potentials}{}
\section{Extensions of Potentials}
\label{section-potentials}
Given a potential $ A \,:\, \Sigma \to \mathbb{R} $, as usual we say that $ \bar{A} \,:\,\hat{\Sigma} \to \mathbb{R} \cup \{ \pm \infty \} $ is an extension of $ A $ to $ \hat{\Sigma} $ whenever $ \bar{A}|_{\Sigma} = A $.
Our aim in this section is to find extension functions on $ \hat{\Sigma} $ which inherit the regularity of $ A $.
Along this section, $ \Sigma $ will always be a general shift space over a countable alphabet.
We will also adopt the name potential for any Borel function $ \bar{A} \,:\,\hat{\Sigma} \to \mathbb{R} \cup \{ -\infty \} $ which is bounded from above.
Although the restriction to $ \Sigma $ of any bounded continuous real-valued function defined on $ \hat{\Sigma} $ is, in fact, a bounded continuous real-valued function on $ \Sigma $, we cannot expect a bounded continuous real-valued function on $ \Sigma $ to always admit a continuous extension on $ \hat{\Sigma} $.
\begin{example}  \label{example-continous-extension}
Let $ \left( u, B_r, B_r, \ldots \right) $ be a point of~$ \partial \Sigma $ and consider a countable partition of $ B_r = \bigsqcup_{k \geq 1 } C_k $ into infinite subsets. 
Let $ \left\{ q_k \right\}_{ k \geq 1 } $ be an enumeration of rational numbers of $ ( 0, 1 ] $.
We introduce the locally constant function $ A \,:\, [ u ] \to [0,1] $  defined on the cylinder set $ [ u ] \in \Sigma $ as
\begin{equation*}
A \left( u a x \right) = A \left( u a \right)
= \left\{
\begin{array}{ c l }
0, & \text{if } a \in \mathscr{A} \setminus B_r \\
q_k, & \text{if } a \in C_k
\end{array} \right., 
\end{equation*}
where $ a \in \mathscr{A} $, $ x \in \Sigma $, and $ u a x  \in \Sigma $.
Due to Tietze extension Theorem, there is a bounded continuous $ A \,:\, \Sigma \to [0,1] $ which extends this locally constant function to~$ \Sigma $. 
However, for any open neighborhood $ Z [ u B_r ; S ] $ of $ \left( u, B_r, B_r, \ldots \right) $,
the oscillation of $ A $ on $ Z [ u B_r ; S ] \cap \Sigma $ is always equal to $ 1 $.
We conclude that there is no continuous extension to $ \hat{\Sigma} $.
\end{example}

On the other hand, upper semi-continuity gives us a lot of freedom to obtain extension functions on $ \hat{\Sigma} $ with the same regularity, as highlighted in the proposition below.
This is a property resulting from the topology of $ \hat{\Sigma} $, especially from the fact that, according to item~\ref{item-convegence-partial-sigma} of Lemma~\ref{lemma-convegence-hat-sigma}, any sequence in $ \partial \Sigma $ converging to a point of the form $ \left( v, B_r, B_r, \ldots\right) \in \partial \Sigma $ is eventually constant.
\begin{proposition}  \label{proposition-Aext-usc}
Let $ A \,:\, \Sigma \to \mathbb{R} $ be an upper semi-continuous potential.
Suppose that $ \bar{A} \,:\,\hat{\Sigma} \to \mathbb{R} \cup \{ -\infty \} $ is an extension of $ A $ fulfilling,
for all $ \left( v, B_r, B_r, \ldots \right) \in \partial \Sigma $,
\begin{equation*}
\limsup_{ B_r \ni i \to \infty } \, \sup A|_{[ v i ]} 
\leq \, \bar{A} \left( v, B_r, B_r, \ldots \right) \,
\leq \sup A|_{ [ v ] }.
\end{equation*}
Then, $ \bar{A} $ is an upper semi-continuous potential.
\end{proposition}

\begin{proof}
Obviously, $ \bar{A} $ is bounded from above. 
Concerning its regularity, for any real number~$ R $, we need to argue that $ \bar{A}^{-1} [-\infty, R) $ is an open set of $ \hat{\Sigma} $. 
Our discussion will be divided into two situations, depending on whether a point of this set belongs to $ \partial \Sigma $ or not.

First, consider a point $ \left( v, B_r, B_r, \ldots \right) $ of $ \partial \Sigma $ verifying
$ \limsup_{B_r \ni i \to \infty} \, \sup A|_{[vi]} \leq  \bar{A} \left( v B_r, B_r, \ldots \right) < R $.
From the definition of superior limit, there exists a finite subset $ S \subset B_r $ such that 
$ \sup A|_{ [ v j ] } < R $, for all $ j \in B_r \setminus S $.
Take the neighborhood $ Z [ v B_r; S ] $ of $ \left( v, B_r, B_r, \ldots \right) $, whose points are of the form 
\begin{equation*}
\left( v, B_r, B_r, \ldots \right),
\qquad 
\left( v, j, w, B_r, B_r, \ldots \right)
\quad \text{or} \quad 
\left( v, j, y_0, y_1, \ldots \right),
\end{equation*}
where $ j \in B_r \setminus S $ is a symbol, $ w \in \mathscr{L} $ is a finite (or empty) word, and $ \left( y_0, y_1, \ldots \right) \in \Sigma $ is an infinite sequence.
In particular, we obtain
\begin{align*}
\bar{A} \left( v, j, w, B_r, B_r, \ldots \right) 
\leq \sup A|_{ [ v j \, w ] }
\leq \sup A|_{ [ v j ] } & < R 
\\
\text{or}\qquad
\bar{A} \left( v, j, y_0, y_1, \ldots \right) 
= A \left( v, j, y_0, y_1, \ldots \right) \leq \sup A|_{ [ v j ] } & < R,
\end{align*}
for every $ j \in B_r \setminus S $.
Hence, $ Z [ v B_r; S ] \subset \bar {A}^{-1} [-\infty, R) $.

Now let $ x \in \Sigma $ be a point such that $ A ( x ) = \bar{A} ( x ) < R $.
Since $ A $ is upper semi-continuous on $ \Sigma $, given $ \varepsilon > 0 $, there exists a neighborhood $ [ x_0, \ldots, x_{n} ] $ of $ x $ such that $ A (z) < R - \varepsilon $ for all $ z \in [ x_0, \ldots, x_{n} ] $.
Consider the neighborhood $ Z [ x_0, \ldots, x_{n} ] $ of $ x $ in~$ \hat{\Sigma} $, whose points are
\begin{equation*}
\left( x_0, \ldots, x_{n}, w, B_r, B_r, \ldots \right)
\qquad \text{or} \qquad 
\left( x_0, \ldots, x_{n}, y_0, y_1, \ldots \right),
\end{equation*}
where $ w \in \mathscr{L} $ is a finite (or empty) word and $ \left( y_0, y_1, \ldots \right) \in \Sigma $ is an infinite sequence.
Note that
\begin{align*}
\bar{A} \left( x_0, \ldots, x_{n}, w, B_r, B_r, \ldots \right) 
\leq \sup A|_{ [ x_0, \ldots, x_{n}, w ] } 
\leq \sup A|_{ [ x_0, \ldots, x_{n} ] } 
\leq R - \varepsilon 
& < R 
\\
\text{or} \quad 
\bar{A} \left( x_0, \ldots, x_{n}, y_0, y_1, \ldots \right) 
= A \left( x_0, \ldots, x_{n}, y_0, y_1, \ldots \right) 
\leq \sup A|_{ [ x_0, \ldots, x_{n} ] } 
\leq R - \varepsilon
& < R.
\end{align*}
Therefore, $ Z [ x_0, \ldots, x_{n} ] \subset \bar {A}^{-1} [-\infty, R) $.
\end{proof}

We will take into account a natural choice of extension potential, which will play a key role in the prove of the \hyperlink{EThm}{Existence Theorem} in Section~\ref{section-existence-maximizing-measure}.
\begin{corollary}  \label{corollary-minimal-Aext}
Let $ A \,:\, \Sigma \to \mathbb{R} \cup \{ -\infty \} $ be an upper semi-continuous potential.
Then, the potential $ \hat{A} \,:\, \hat{\Sigma} \longrightarrow \mathbb{R} \cup \{ -\infty \} $ given as
\begin{align*}
x \in \Sigma & \longmapsto \hat{A} (x) = A (x) 
\qquad \text{and}\\ 
\left( v, B_r, B_r, \ldots \right) \in \partial \Sigma 
& \longmapsto \hat{A} \left( v, B_r, B_r, \ldots \right) = \limsup_{B_r \ni i \to \infty} \, \sup A|_{[ v i ]}
\end{align*}
is the minimal upper semi-continuous extension of $ A $ to $ \hat{\Sigma} $.
\end{corollary}

\begin{proof}
Due to Proposition~\ref{proposition-Aext-usc}, we already know that $ \hat{A} $ is an upper semi-continuous potential.
Let $ \bar{A} \,:\, \Sigma \to \mathbb{R} \cup \{ -\infty \} $ be any upper semi-continuous extension of~$ A $.
We claim that $ \hat{A}|_{ \partial \Sigma } \leq \bar{A}|_{ \partial \Sigma } $.
As matter of fact, given any real number $ R $ such that $ \bar{A} \left( v, B_r, B_r, \ldots \right) < R $, there exists a finite subset $ S_R \subset B_r $ that verifies $ Z [ v B_r ; S_R ] \subset \bar{A}^{-1} [-\infty, R) $.
Since $ [ v i ] \subset Z [ v B_r ; S_R ] $ for $ i \in B_r \setminus S_R $, note that  
$ \sup A|_{ [ v i ] } 
\leq \sup \bar{A}|_{ Z [ v B_r ; S_R ] }
\leq R $ 
whenever $ i \in B_r \setminus S_R $. 
By definition, we have
$ \hat{A} \left( v, B_r, B_r, \ldots \right) 
= \limsup_{ B_r \ni i \to \infty } \sup A|_{ [ v i ] } 
\leq R $.
Hence, $ \hat{A} \left( v, B_r, B_r, \ldots \right) $ became a lower bound for any real number $ R > \bar{A} \left( v, B_r, B_r, \ldots \right) $.
We conclude that $ \hat{A} \left( v, B_r, B_r, \ldots \right) \leq \bar{A} \left( v, B_r, B_r, \ldots \right) $.
\end{proof}

A well studied class of potentials in the literature \cite{JMU:ETDS06, BG:BBMS10, BF:ETDS13, GG:arXiv24} is the coercive one, namely, those potentials $ A $ fulfilling
$ \lim_{i \to \infty} \, \sup A|_{[i]} = - \infty $.
The corresponding minimal extension $ \hat{A} $ exhibits in the coercive framework a convenient behavior as described in the following example.
\begin{example}
Assume that $ A $ is a coercive upper semi-continuous potential and consider the minimal upper semi-continuous extension $ \hat{A} $.
It is easy to see that 
\begin{equation*}
\int \hat{A} \, d \delta_{ \left( B_r, B_r, \ldots \right) }
= \hat{A} \left( B_r, B_r, \ldots \right) 
= \limsup_{ B_r \ni i \to \infty } \, \sup A|_{[i]} 
= \lim_{ i \to \infty } \sup A|_{[i]} 
= -\infty
\end{equation*}
for every $ \left( B_r, B_r, \ldots \right) \in \hat{\mathcal{L}}_0 $.
In particular, the coercive assumption prevents the invariant probabilities giving mass $ 1 $ to $ \partial \Sigma $ from being maximizing measures.
\end{example}

In Section~\ref{section-existence-maximizing-measure}, we will show that such a behavior may be observed with respect to more general potentials $ A $: those satisfying
$ \limsup_{ i \to \infty } \, \sup A|_{[i]} < \beta (A) $,
a key assumption in \hyperlink{EThm}{Existence Theorem}.
In order to do that, we will study in \hyperlink{sec-ergodic-optimization}{next section} the main aspects of ergodic optimization theory on blur shifts.

\hypertarget{sec-ergodic-optimization}{}
\section{Ergodic Optimization}
\label{section-ergodic-optimization}

Throughout this section, let $ \Sigma $ be any shift over countable alphabet, 
assume that $ A \,:\, \Sigma \to \mathbb{R} $ is a potential and 
consider $ \bar{A} \,:\, \hat{\Sigma} \to \mathbb{R} \cup \{ - \infty \} $ any extension potential of~$ A $.
We introduce the following ergodic maximizing constants
\begin{align*}
\beta \left( A \right)
& = \sup_{ \mu \in \mathcal{M} \left( \Sigma, \sigma \right) } \,
\int_{ \Sigma } A \, d \mu, 
\\[0.1cm]
\hat{ \beta } \big( \bar{A} \big) 
& = \sup_{ \hat{\mu} \in \mathcal{M} \big( \hat{\Sigma}, \hat{\sigma} \big) } \, 
\int_{ \hat{\Sigma} } \bar{A} \, d \hat{\mu},
\\[0.1cm]
\text{and} \quad
\max \, \bar{A}|_{\hat{\mathcal{L}}_{0} } 
& = \max_{ 1 \leq r \leq s } \, 
\int \bar{A} \, d \delta_{ \left( B_r, B_r, \ldots \right) }.
\end{align*}
As a notational shorthand, we may also use $ \vee $ for the maximum between values.
The above constants are related as follows.
\begin{lemma}  \label{lemma-ergodic-maximizing-value}
\enspace
$ \displaystyle
\hat{\beta} \big( \bar{A} \big) 
\, = \,
\beta \left( A \right)
\, \vee \,
\max \, \bar{A}|_{ \hat{\mathcal{L}}_{0} }.
$
\end{lemma}

\begin{proof}
Remember that
$ \mathcal{M} \left( \Sigma, \sigma \right) \subset \mathcal{M} \big( \hat{\Sigma}, \hat{\sigma} \big) $.
Besides, 
$ \bar{A} \left( B_r, B_r, \ldots \right) = \int \bar{A} \, d \delta_{ \left( B_r, B_r, \ldots \right) } $ 
and $ \delta_{ \left( B_r, B_r, \ldots \right) } $ is an $\hat{\sigma} $-invariant probability measure, for every $ 1 \leq r \leq s $.
It immediate from these facts that
$ \beta \left( A \right) \vee \max \, \bar{A}|_{\hat{\mathcal{L}}_{0} } \leq \hat{\beta} \big( \bar{A} \big) $.
On the other hand, the convex characterization of Proposition~\ref{proposition-decomposition-blur-invariant-measures} implies
\begin{align*}
\hat{\beta} \big( \bar{A} \big) 
& = \sup \, \bigg\{ \,
t \int_{ \Sigma } A \, d \dot{\mu} \,+\,
( 1 - t ) \int_{ \partial \Sigma } \bar{A} \, d \hat{\nu}
\,:\, \\
& \hspace*{1.2cm}
t \in [0,1], \
\dot{\mu} \in \mathcal{M} \left( \Sigma, \sigma \right) 
\text{ and }
\hat{\nu} \in \mathrm{Conv} \left( \, \left\{ \delta_{ \left( B_r, B_r, \ldots \right) } \,:\, 1 \leq r \leq s  \right\} \, \right)
\, \bigg\} 
\\
& \leq \beta \left( A \right) \, \vee \,
\max \, \Bigg\{ \,
\sum_{ r = 1 }^{ s } \alpha_r \, \bar{A} \left( B_r, B_r, \ldots \right) \,:\,
\sum_{ r = 1 }^{ s } \alpha_r = 1 
\, \Bigg\}.
\end{align*}
Clearly, the last term is less than or equal to
$ \max \, \bar{A}|_{ \hat{\mathcal{L}}_{0}  }$.
\end{proof}

Henceforward, we take into account the maximizing sets
\begin{align*}
\mathcal{M}_{\max} \left( A \right)
& = \left\{ \mu \in \mathcal{M} \left( \Sigma, \sigma \right) \,:\, \int_{ \Sigma } A \, d \mu = \beta ( A ) \right\},
\\[0.1cm]
\hat{\mathcal{M}}_{\max} \big( \bar{A} \big) 
& = \left\{ \hat{\mu} \in \mathcal{M} \big( \hat{\Sigma}, \hat{\sigma} \big) \,:\, \int_{ \hat{\Sigma} } \bar{A} \, d \hat{\mu} = \hat{\beta} \big( \bar{A} \big) \right\}, 
\quad \text{and}
\\[0.1cm]
\hat{\mathcal{D}}_{\max} \big( \bar{A} \big) 
& = \left\{ \delta_{ \left( B_r, B_r, \ldots \right) } \,:\, 
\bar{A} \left( B_r, B_r, \ldots \right) = \max \, \bar{A}|_{\hat{\mathcal{L}}_{0} } \right\}.
\end{align*}
It is clear that $ \hat{\mathcal{D}}_{\max} $ is a non-empty set.
Nevertheless, the convex sets $ \mathcal{M}_{\max} \left( A \right) $ and $ \hat{\mathcal{M}}_{\max} \big( \bar{A} \big) $ may be empty or not.
The maximizing sets are associated according to the following result.
\begin{proposition}  \label{proposition-behavior-set-maximizing-measure} \ 
\begin{enumerate}[label=(\roman*)]
\item\label{item-beta-greater-max-Br}
If $ \max \, \bar{A}|_{ \hat{\mathcal{L}}_{0} } < \beta \left( A \right) $,
then 
$ \hat{\mathcal{M}}_{\max} \big( \bar{A} \big) = \mathcal{M}_{\max} \left(  A \right) $.
\\[-0.2cm]
\item\label{item-beta-equal-max-Br} 
If $ \max \, \bar{A}|_{ \hat{\mathcal{L}}_{0} } = \beta \left( A \right) $,
then
$ \hat{\mathcal{M}}_{\max} \big( \bar{A} \big) =
\mathrm{Conv} \big( \, \mathcal{M}_{\max} \left( A \right) \, \sqcup \, \hat{\mathcal{D}}_{\max} \big( \bar{A} \big) \, \big) $.
\\[-0.2cm]
\item\label{item-beta-less-max-Br} 
If $ \max \, \bar{A}|_{ \hat{\mathcal{L}}_{0} } > \beta \left( A \right) $,
then
$ \hat{\mathcal{M}}_{\max} \big( \bar{A} \big) = \mathrm{Conv} \big( \, \hat{\mathcal{D}}_{\max} \big( \bar{A} \big) \, \big) $.
\end{enumerate}
\end{proposition}

\begin{proof}
First note that Lemma~\ref{lemma-ergodic-maximizing-value} and Proposition~\ref{proposition-decomposition-blur-invariant-measures} guarantee that any maximizing measure $ \hat{\mu} \in \hat{\mathcal{M}}_{\max} \big( \bar{A} \big) $ (when it exists) satisfies
\begin{equation}  \label{equation-integral-decomposition-maximizing-measure}
\beta \left( A \right) \, \vee \, \max \, \bar{A}|_{ \hat{\mathcal{L}}_{0} } 
= \int_{ \hat{\Sigma} } \bar{A} \, d \hat{\mu} 
= t \int_{\Sigma} A \, d \dot{\mu} + (1-t) \, \sum_{ r = 1 }^{ s } \alpha_r \, \bar{A} \left( B_r, B_r, \ldots \right),
\end{equation}
where 
$ t \, \dot{\mu} + (1-t) \, \sum_{r=1}^{s} \alpha_r \, \delta_{ \left( B_r, B_r, \ldots \right) } $, with 
$ t \in [0,1] $, 
$ \dot{\mu} \in \mathcal{M} \left( \Sigma, \sigma \right) $,
$ \alpha_r \geq 0 $ and $ \sum_{r} \alpha_r = 1 $, 
is the convex representation of $ \hat{\mu} $.

\smallskip
{\itshape Item~\ref{item-beta-greater-max-Br}.}
From Lemma~\ref{lemma-ergodic-maximizing-value}, we have $ \hat{\beta} \big( \bar{A} \big) = \beta \left( A \right) $.
It is straightforward that
$ \mathcal{M}_{\max} \left( A \right) \subset \hat{\mathcal{M}}_{\max} \big( \bar{A} \big) $.
On the other hand, assume that there exists a maximizing probability $ \hat{\mu} \in \hat{\mathcal{M}}_{\max} \big( \bar{A} \big) $
(otherwise, the result is trivial).
By~\eqref{equation-integral-decomposition-maximizing-measure}, it follows that
\begin{equation*}
\beta \left( A \right) 
= t \int_{\Sigma} A \, d \dot{\mu} + (1-t)\, \sum_{ r = 1 }^{ s } \alpha_r \, \bar{A} \left( B_r, B_r, \ldots \right) 
\leq t \, \beta \left( A \right) + (1-t) \, \max \, \bar{A}|_{ \hat{\mathcal{L}}_{0} }.
\end{equation*}
Due to the hypothesis, this inequality only holds when $ t = 1 $.  
Thus, the maximizing measure $ \hat{\mu} = \dot{\mu} \in \mathcal{M} \left( \Sigma, \sigma \right) $ verifies $ \int A \, d \dot{\mu} =  \beta \left( A \right) $,
which means that $ \hat{\mathcal{M}}_{\max} \big( \bar{A} \big) \subset \mathcal{M}_{\max} \left( A \right) $.

\smallskip
{\itshape Item~\ref{item-beta-equal-max-Br}.}
From Lemma~\ref{lemma-ergodic-maximizing-value}, 
$ \hat{\beta} \big( \bar{A} \big) 
= \beta \left( A \right)  
= \max \, \bar{A}|_{ \hat{\mathcal{L}}_{0} } $.
These equalities implies that $ \hat{\mathcal{M}}_{\max} \big( \bar{A} \big) $ contains both $ \mathcal{M}_{\max} \left( A \right) $ and $ \hat{\mathcal{D}}_{\max} \big( \bar{A} \big) $.
In particular, $ \hat{\mathcal{M}}_{\max} \big( \bar{A} \big) $ is non-empty and, because of its convexity, 
$ \mathrm{Conv} \big( \, \mathcal{M}_{\max} \left( A \right) \, \sqcup \, \hat{\mathcal{D}}_{\max} \big( \bar{A} \big) \, \big) 
\subset \hat{\mathcal{M}}_{\max} \big( \bar{A} \big) $.
Given $ \hat{\mu} \in \hat{\mathcal{M}}_{\max} \big( \bar{A} \big) $,
since $ \hat{\mu} = t \, \dot{\mu} + (1-t) \, \sum_{r=1}^{s} \alpha_r \, \delta_{ \left( B_r, B_r, \ldots \right) } $, 
it remains to show that  
$ \dot{\mu} \in \mathcal{M}_{\max} \left( A \right) $ and 
$ \sum_{ r = 1 }^{ s } \alpha_r \, \delta_{ \left( B_r, B_r, \ldots \right) } \in \mathrm{Conv} \big( \, \hat{\mathcal{D}}_{\max} \big( \bar{A} \big) \, \big) $.
If either $ t = 0 $ or $ t = 1 $, this is obvious.
For $ t \in ( 0, 1 ) $, it follows from~\eqref{equation-integral-decomposition-maximizing-measure} that
\begin{align*}
\int_{\Sigma} A \, d \dot{\mu} 
= \dfrac{1}{t} \left[ \,
\beta \left( A \right) 
- (1-t) \, \sum_{ r = 1 }^{ s } \alpha_r \, \bar{A} \left( B_r, B_r, \ldots \right) 
\, \right]
& \geq \beta (A) 
\\
\text{and} \quad
\sum_{ r = 1 }^{ s } \alpha_r \, \bar{A} \left( B_r, B_r, \ldots \right)
= \dfrac{1}{1-t} \left[ \,
\max \, \bar{A}|_{ \hat{\mathcal{L}}_{0} } 
- t \int_{\Sigma} A \, d \dot{\mu} 
\, \right]
& \geq \max \, \bar{A}|_{ \hat{\mathcal{L}}_{0} }.
\end{align*}
It is easy to see that the above inequalities are in fact equalities, from which we immediately conclude the discussion.

\smallskip
{\itshape Item~\ref{item-beta-less-max-Br}.}
Lemma~\ref{lemma-ergodic-maximizing-value} imply $ \hat{\beta} \big( \bar{A} \big) = \max \, \bar{A}|_{ \hat{\mathcal{L}}_{0} } $.
Consequently, $ \hat{\mathcal{D}}_{\max} \big( \bar{A} \big) \subset \hat{\mathcal{M}}_{\max} \big( \bar{A} \big) $
and $ \hat{\mathcal{M}}_{\max} \big( \bar{A} \big) $ is non-empty.
By convexity,
$ \mathrm{Conv} \big( \, \hat{\mathcal{D}}_{\max} \big( \bar{A} \big) \, \big) \subset \hat{\mathcal{M}}_{\max} \big( \bar{A} \big) $. 
Given any $ \hat{\mu} \in \hat{\mathcal{M}}_{\max} \big( \bar{A} \big) $, 
it follows from~\eqref{equation-integral-decomposition-maximizing-measure} that
\begin{equation*}
\max \, \bar{A}|_{ \hat{\mathcal{L}}_{0} }
= t \int_{\Sigma} A \, d \dot{\mu} + (1-t) \, \sum_{ r = 1 }^{ s } \alpha_r \, \bar{A} \left( B_r, B_r, \ldots \right) 
\leq t \, \beta \left( A \right) + (1-t) \, \max \, \bar{A}|_{ \hat{\mathcal{L}}_{0} },
\end{equation*} 
which, thanks to the hypothesis, only holds if $ t = 0 $.
Hence,
$ \bar{A} \left( B_r, B_r, \ldots \right) = \max \, \bar{A}|_{ \hat{\mathcal{L}}_{0} } $ whenever~$ \alpha_r > 0 $, 
and we obtain that 
$ \hat{\mu} = \sum_{ r = 1 }^{ s } \alpha_r \, \delta_{ \left( B_r, B_r, \ldots \right) } \in \mathrm{Conv} \big( \, \hat{\mathcal{D}}_{\max} \big( \bar{A} \big) \, \big) $.
\end{proof}

\begin{remark}
Since the set $ \hat{\mathcal{D}}_{\max} $ is always non-empty, we can conclude that the converse of item~\ref{item-beta-greater-max-Br} of Proposition~\ref{proposition-behavior-set-maximizing-measure} is also true.
\end{remark}

\section{Existence of Maximizing Measure}
\label{section-existence-maximizing-measure}
This section is dedicated to prove the \hyperlink{EThm}{Existence Theorem} via blur shift compactification method.

\begin{proof}[Proof of the \hyperlink{EThm}{Existence Theorem}]
Due to Corollary~\ref{corollary-blur-measure-compact}, $ \mathcal{M} \big( \hat{\Sigma}, \hat{\sigma} \big) $ is a compact set.
Note that the upper semi-continuous extension potential $ \hat{A} $ (guaranteed by Corollary~\ref{corollary-minimal-Aext}) induces an upper semi-continuous application 
$ \hat{\mu} \in \mathcal{M} \big( \hat{\Sigma}, \hat{\sigma} \big) 
\longmapsto \int_{ \hat{\Sigma} } \hat{A} \, d \hat{\mu} $.
Therefore, $ \hat{\mathcal{M}}_{\max} \big( \hat{A} \big) $ is non-empty.
From the key assumption on the potential, we obtain
\begin{equation*}
\max \, \hat{A}|_{ \hat{\mathcal{L}}_{0} } 
= \max \left\{ \, \limsup_{B_r \ni i \to \infty} \, \sup A|_{[i]} \,:\, 1 \leq r \leq s \, \right\} 
= \limsup_{i \to \infty} \, \sup A|_{[i]} 
< \beta (A).
\end{equation*}
Hence, item~\ref{item-beta-greater-max-Br} of Proposition~\ref{proposition-behavior-set-maximizing-measure} implies that
$ \mathcal{M}_{\max} \left(  A \right) = \hat{\mathcal{M}}_{\max} \big( \hat{A} \big) $ is non-empty.
The last statement of the theorem will follow from Proposition~\ref{proposition-finite-shift-support}. 
\end{proof}

\begin{proposition}  \label{proposition-finite-shift-support}
Let $ A \,:\, \Sigma \to \mathbb{R} $ be an upper semi-continuous potential that admits a maximizing probability $ \mu \in \mathcal{M} \left( \Sigma, \sigma \right) $ and fulfills both conditions 
\begin{equation*}
\sup A \leq  \beta \big( A \big)
\qquad \text{and} \qquad
\limsup_{ i \to \infty } \, \sup A|_{[i]} < \beta ( A ),
\end{equation*}
Then, there exists a common finite alphabet subshift of~$ \Sigma $ that contains the support of any maximizing measure for $ A $. 
\end{proposition}

\begin{proof}
\let\qed\relax
Since $ \beta (A) \leq \sup A $ in general, the hypothesis means that $ \beta (A) = \sup A $.
First, we claim that the support of any maximizing probability for $ A $ lies in $ A^{-1} \left( \, \beta (A) \, \right) $.
In fact, for $ \mu \in \mathcal{M}_{\max} \left( A \right) $, we clearly have
\begin{equation*}
\int_{\Sigma} A \, d\mu = \beta (A) =  \sup A
\qquad \Longrightarrow \qquad
\mu \big( A^{-1} \left( \, \beta (A) \, \right) \big) = 1,
\end{equation*}
By the upper semi-continuity of~$ A $, it is immediate that $ A^{-1} \left( \, \beta (A) \, \right) $ is a closed subset of $ \Sigma $, which shows our claim.

Let us prove the desired result by contrapositive argument. 
Assume that there exists an infinite sub-alphabet $ \mathscr{C} $ of $ \mathscr{L}_1 $ 
such that
each cylinder set $ [i] $, for $ i \in \mathscr{C} $, has non-empty intersection with the support of some maximizing measure for~$ A $. 
By the previous claim, 
$ [i] \, \cap \, {A}^{-1} \left( \, \beta (A) \, \right) $ is also non-empty, so that
$ \sup A|_{[i]} \geq \beta (A) $ for every $ i \in \mathscr{C} $.
We conclude that  
\begin{equation*}
\limsup_{ i \to \infty } \, \sup A|_{[i]}
\geq \limsup_{ \mathscr{C} \ni i \to \infty } \, \sup A|_{[i]}
\geq \beta ( A ).
\tag*{\qedsymbol}
\end{equation*}
\end{proof}

We end this section by showing that the converse of the above proposition is false; more precisely, it may happen that the supports of all maximizing measures lie in a certain finite alphabet subshift and yet $ \limsup_{ i \to \infty } \, \sup A|_{[i]}
= \beta ( A ) $.
\begin{example} \label{example-converse-finite-shift-support}
Let $ y \in \Sigma $ be a periodic point, with period $ p $, and consider the associated periodic probability
$ \mu_{y} = \frac{1}{ p } \, \sum_{ i = 0 }^{ p - 1 } \delta_{ \sigma^i ( y ) } $.
Introduce the bounded continuous potential $ A \,:\, \Sigma \to \mathbb{R} $ given as
\begin{equation*}
A (x) = A \left( x_0, x_1, \ldots \right) = - \dfrac{1}{x_0} \, d \left( \mathrm{orb} \, (y), x \right),
\end{equation*}
where $ d $ denotes the usual distance on $ \Sigma $, namely, two infinite sequences of $ \Sigma $ are far from each other~$ \nicefrac{1}{2} $ raised to the power of their first disagreement.
Obviously, $ \beta (A) = 0 $ and $ \mathcal{M}_{\max} \left( A \right) = \left\{ \mu_{y} \right\} $.
In particular, $ \sup A \leq \beta \big( A \big) $ and the subshift over a finite alphabet $ \mathrm{orb} \, (y) $ is the support of the unique maximizing probability for $ A $. 
On the other hand, 
\begin{align*}
\limsup_{ i \to \infty } \, \sup A|_{[i]}
= - \lim_{ i \to \infty } \, \dfrac{1}{i} \, \inf_{x \in [i]} \,  d \left( \mathrm{orb} \, (y), \, x \right)
= 0
= \beta (A).
\end{align*}
Note that, in this case, the minimal extension $ \hat{A} $ fulfills $ \hat{A}|_{ \hat{\mathcal{L}}_{0} } \equiv \beta ( A ) $ and, 
by item~\ref{item-beta-equal-max-Br} of Proposition~\ref{proposition-behavior-set-maximizing-measure}, 
$ \hat{\mathcal{M}}_{\max} \big( \hat{A} \big) = 
\mathrm{Conv} \left( \, \left\{ \mu_{y} \right\} \, \sqcup \, \left\{ \delta_{ \left( B_r, B_r, \ldots \right) } \,:\, 1 \leq r \leq s  \right\}  \, \right) $.
\end{example}

\hypertarget{sec-non-markov-example}{}
\section{Non-Markovian Example}
\label{section-non-markov-example}

We exhibit in this section a countable alphabet non-Markovian shift to which we may apply the results of this paper.
Recall that
a countable alphabet shift is called Markovian if there exists countable zero-one matrix $ M $ such that any point $ \left( x_0, x_1, \ldots \right) $ obeys $ M (x_i, x_{i+1}) = 1 $ for all $ i \geq 0 $
(for a comprehensive study, see, for instance, \cite{Kit:Spr98}).
A key aspect of our analysis will be the specification property, which plays a crucial role in establishing the denseness of periodic measures. Observe that a shift space has the specification property if, given any finite collection of cylinder sets and sufficiently large gaps between them, one can construct a point that visits each cylinder set at the prescribed times while filling the gaps with arbitrary transitional words.

\begin{definition} \label{definition-non-markov-example}
Consider a partition $ \mathbb{Z}_{+} = \bigsqcup_{k=1}^{\infty} Z_{k} $ such that 
$ Z_1 = \{ 1 \} $,
$ Z_2 = \{ 2, 3 \} $  and 
$ Z_k $ are finite sets with increasing cardinality, i.e., $ | Z_k | < | Z_{k+1} | $ for all $ k \geq 1 $.
Let $ \Lambda $ be a (one-sided) shift space defined over a countable alphabet $ \mathbb{Z}_{+} $ that satisfies the following conditions:
\begin{enumerate}[label=(C\arabic*)]
\item\label{item-C1} 
For every point $ x \in \Lambda $, there exists an integer $ N_x \geq 0 $, such that $ \sigma^n (x) \notin [1]$ for any $ n \geq N_x $.
In other terms, the symbol $1$ appears a finite number of times for each sequence.
\item\label{item-C2} 
For $ a \in Z_k $ and $ k > 1 $, the predecessor set $ \mathscr{P} (a) $ is included in $ Z_k $, with the exception of $ a = \min \, Z_{k} $ for which one also admits $ 1 $ as an extra predecessor.
In particular, the follower set $ \mathscr{F}_1 (1) = \{ \min \, Z_k \,:\, k > 1 \} $.
\item\label{item-C3} 
For $ b \in Z_k $ and $ k > 1 $, the follower set $ \mathscr{F}_1 (b) $ is included in $ Z_k $, with the exception of $ b = \max \, Z_{k} $ for which one also admits $ 1 $ as an extra follower.
In particular, the predecessor set $ \mathscr{P} (1) = \{\max \, Z_{k} \,:\, k > 1 \} $.
\item\label{item-C4} 
The subshift $ \Lambda \cap Z_{2}^{ \mathbb{N} } $ is the sofic even shift and,
for any $ k > 2 $, each finite alphabet subshift $ \Lambda \cap Z_{k}^{ \mathbb{N} } $ verifies the specification property.
In particular, the restriction $ \sigma|_{ \Lambda \cap Z_{k}^{ \mathbb{N} } } $ is transitive, for all $ k > 1 $.
\item\label{item-C5} 
Any word of the language of $ \Lambda $ with first symbol $ \min \, Z_{k} $ and last symbol $ \max \, Z_{k} $ must have length at least $|Z_{k}|$.
\end{enumerate}
\end{definition}

\begin{figure}[ht]
\centering
\begin{tikzpicture}[scale=0.8, every node/.style={transform shape}]
\tikzstyle{vertice} = [rectangle, fill=white, line width=0.8pt, draw=black!85, rounded corners=6pt, inner sep=5pt, font=\normalsize]
\tikzstyle{particao} = [font=\small]
\tikzstyle{moldura} = [fill=gray!70, rounded corners, font=\normalsize]
\tikzstyle{aresta} = [-stealth, draw=black!90, line width=0.8pt]

\node[vertice, circle] (1) at (0,0) {$ 1 $};
\node[particao, left = 1pt of 1.west] (Z1) {$ Z_{1} $};
\begin{scope}[on background layer]
\node[moldura, fit = {(1) (Z1)}] {}; 
\end{scope}

\node[vertice] (2) at (3,0) {$ \min \, Z_{2} $};
\node[below = 0pt of 2, inner sep=0pt] (dotsZ2) {$ \vdots $};
\node[vertice, below = 20pt of 2] (3) {$ \max \, Z_{2} $};
\node[particao, above = 1pt of 2.north] (Z2) {$ Z_2 $};
\begin{scope}[on background layer]
\node[moldura, fit = { (Z2) (2) (dotsZ2) (3) }] {}; 
\end{scope}

\node[vertice] (4) at (6,1.5) {$ \min \, Z_{3} $};
\node[below = 10pt of 4] (dotsZ3) {$ \vdots $};
\node[vertice, below = 50pt of 4] (5) {$ \max \, Z_{3} $};
\node[particao, above = 1pt of 4] (Z3) {$ Z_3 $};
\begin{scope}[on background layer]
\node[moldura, fit = { (Z3) (4) (dotsZ3) (5) }] {}; 
\end{scope}

\node[vertice] (6) at (9,3) {$ \min \, Z_{4} $};
\node[below = 25pt of 6] (dotsZ4) {$ \vdots $};
\node[vertice, below = 80pt of 6] (7) {$ \max \, Z_4 $};
\node[particao, above = 1pt of 6] (Z4) {$ Z_4 $};
\begin{scope}[on background layer]
\node[moldura, fit = { (Z4) (6) (dotsZ4) (7) }] {}; 
\end{scope}

\node[] (101) at (12,4.5) {$ \ldots $};
\node[] (102) at (12,-3) {$ \ldots $};

\draw[aresta] (1) to[out = 0, in = 180] (2);
\draw[aresta, rounded corners = 5pt] (3) -| (1.-24);

\draw[aresta, rounded corners = 10pt] (1.45) |- (4);
\draw[aresta, rounded corners = 10pt] (5) |- (2,-2.2) -| (1.-68);

\draw[aresta, rounded corners = 10pt] (1.90) |- (6);
\draw[aresta, rounded corners = 10pt] (7) |- (2,-2.6) -| (1.-112);

\draw[aresta, rounded corners = 10pt] (1.135) |- (101);
\draw[aresta, rounded corners = 10pt] (102) -| (1.-156);
\end{tikzpicture}
\caption[Graph representation transition of $ \Lambda $]
{Graph representation of conditions \ref{item-C2} and \ref{item-C3} in Definition~\ref{definition-non-markov-example}}
\label{figure-non-markov-example}
\end{figure}
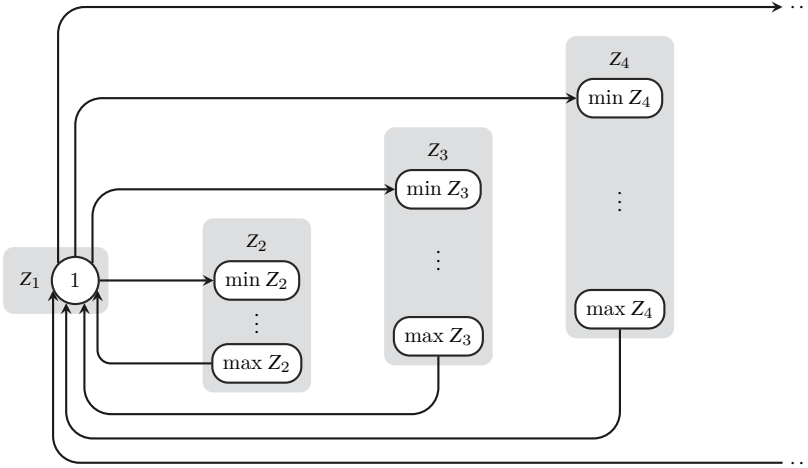

We can explicitly exhibit the predecessor and follower relations of \ref{item-C2} and \ref{item-C3} as arrow connections on a graph given in Figure~\ref{figure-non-markov-example}. 
Moreover, conditions~\ref{item-C4} and~\ref{item-C5} may be interpreted as the subshift $ \Lambda \cap Z_{k}^{\mathbb{N}} $ verifying certain regularity properties that ensure well-behaved dynamics within each partition class while maintaining global non-Markovian structure. Specifically, condition~\ref{item-C4} establishes that the restriction to each $Z_k$ with $k \geq 2$ inherits the specification property from finite alphabet theory, thereby guaranteeing transitivity and enabling the construction of connecting words between arbitrary cylinder sets.
Meanwhile, condition~\ref{item-C5} imposes a minimal length constraint that prevents trivial loops and ensures that transitions between the extremal elements of each partition class cannot occur too rapidly, thus preserving the hierarchical structure encoded in the increasing cardinalities $|Z_k|$.
(Observe that, if we were to impose only condition~\ref{item-C5} to define all forbidden words, we would be led to consider a subshift of finite type that clearly satisfies the specification property, which thus constitutes a natural example of finite alphabet shift to be considered at any level.)
Together, these conditions create a delicate balance: while the global system exhibits non-Markovian behavior, the local dynamics within each class retain sufficient regularity to support the specification property and its associated consequences for the ergodic theory of the system.

In the following results, we will show that any shift $ \left( \Lambda, \sigma \right) $ verifying the conditions in the above Definition~\ref{definition-non-markov-example} possesses the properties:
\begin{itemize}
\item it is non-Markovian
\enspace (Lemma~\ref{lemma-non-markov});
\item it obeys the \hyperlink{fcpa}{finite cyclic predecessor assumption}
\enspace (Lemma~\ref{lemma-non-markov-fcpa});
\item it is non-locally compact and the set $ \partial \Lambda \setminus \hat{\mathcal{L}}_{0} $ is non-empty 
\enspace (Lemma~\ref{lemma-non-markov-non-loc-compact});
\item it is topologically transitive
\enspace (Lemma~\ref{lemma-non-markov-topological-transitive});
\item it satisfies the \hyperlink{dpm}{denseness of periodic measures}
\enspace (Corollary~\ref{corollary-non-markov-dpm}).
\end{itemize}

\begin{lemma}  \label{lemma-non-markov}
The shift space $ \Lambda $ is non-Markovian.
\end{lemma}

\begin{proof}
Since $ \Lambda $ contains the sofic even shift (property~\ref{item-C4}), it is non-Markovian.
\end{proof}

\begin{lemma}  \label{lemma-non-markov-fcpa}
For the shift space $ \Lambda $, we have
\begin{enumerate}[label=(\roman*)]
\item \label{item-fcpa-1}
$ \mathscr{P} (1) \cap \mathscr{F}_{m} (1) $ is finite for all $ m \geq 1 $;
\item \label{item-fcpa-n}
$ \mathscr{P} (a) $ is finite for every $ a \geq 2 $.
\end{enumerate}
In particular, $ \Lambda $ obeys the \hyperlink{fcpa}{finite cyclic predecessor assumption}.
\end{lemma}

\begin{proof}
Item~\ref{item-fcpa-1} is a direct consequence of properties~\ref{item-C2}, \ref{item-C3} and~\ref{item-C5}, and item~\ref{item-fcpa-n} is a reinterpretation of \ref{item-C2}.
\end{proof}

\begin{lemma}  \label{lemma-non-markov-non-loc-compact}
The shift space $ \Lambda $ is non-locally compact.
Moreover, for any associated blur shift $ \big( \hat{\Lambda}, \hat{\sigma} \big) $ with finite resolution, the set $ \bigsqcup_{ k \geq 1 } \hat{\mathcal{L}}_k $ is non-empty.
\end{lemma}

\begin{proof}
Consider any cylinder set $ \left[ x_0, \ldots, x_{n} \right] $.
Due to \ref{item-C4} and \ref{item-C3}, there exists a word $ w $ of the language of $ \Lambda $ with first symbol $ x_{n} $ and last symbol $ 1 $.
We can define a sequence of points
\begin{equation*}
y^{k} := \left( x_0, \ldots, x_{n-1}, \, w, \, \min \, Z_{k}, \ldots \right) 
\in \left[ x_0, \ldots, x_{n} \right]
\subset \Lambda
\qquad 
\text{for any } k > 1,
\end{equation*}
which has no accumulation point on $ \Lambda $.
Therefore, the neighborhood basis consisting of cylinder sets is not (sequentially) compact. 

In order to prove the second assertion, we can view the previous sequence~$ \{ y^{k} \}_{ k > 1 } $ 
as being defined in a blur shift space $ \hat{\Lambda} $ with resolution $ \mathscr{V} = \left\{ B_1, \ldots, B_s \right\} $ associated with $ \Lambda $.
In~particular, there exists some $ B_r $ such that $ \left\{ \min \, Z_{k} \,:\, k > 1 \right\} \cap B_r $ is infinite.
Using item~\ref{item-convegence-partial-sigma} of Lemma~\ref{lemma-convegence-hat-sigma}, there exists a limit point
$ \left( x_0, \ldots, x_{n-1}, \, w, \, B_r, B_r, \ldots \right)
\in \bigsqcup_{ k \geq 1 } \hat{\mathcal{L}}_k $
for the subsequence 
$ y^{k} = \left( x_0, \ldots, x_{n-1}, \, w, \, \min \, Z_{k}, \ldots \right) $,
as $ \min \, Z_{k} \in B_r $ with $ k \to \infty $.
\end{proof}

\begin{remark}
Note that any associated blur shift $ \big( \hat{\Lambda}, \hat{\sigma} \big) $ with finite resolution is non-trivial, in the sense that the set $ \bigsqcup_{ k \geq 1 } \hat{\mathcal{L}}_k $ is non-empty.
Nevertheless, this set has always null measure, as a consequence of item~\ref{item-sum-lk-null-measure} of Proposition~\ref{proposition-equivalent-blur-measure-compact}.
\end{remark}

\begin{lemma}  \label{lemma-non-markov-topological-transitive}
The shift $ \left( \Lambda, \sigma \right) $ is topologically transitive.
\end{lemma}

\begin{proof}
Let $ [ x_0, \ldots, x_{n} ] $ and $ [ y_0, \ldots, y_{m} ] $ be two non-empty cylinder sets in the shift~$ \Lambda $.
From \ref{item-C2}, \ref{item-C3} and \ref{item-C4} it is straightforward to guarantee the existence of the following words: 
\begin{itemize}
\item
$ u $, with first symbol $ x_n $ and last symbol $ 1 $, such that the word $ \left( x_0, \ldots, x_{n-1}, \, u  \right) $ belongs to the language of $ \Lambda $ --
in particular, $ u = 1 $ when $ x_n = 1 $;
\item
$ v $, with first symbol $ 1 $ and last symbol $ y_0 $, such that the word $ \left( v, \, y_1, \ldots, y_{m} \right) $ belongs to the language of $ \Lambda $ -- in particular, $ v = 1 $ when $ y_0 = 1 $.
\end{itemize}
Hence, there exists a point
$ z = \left( x_0, \ldots, x_{n-1}, \, u, \, 2, 3, \, v, \, y_1, \ldots, y_{m}, \ldots \right) 
\in [ x_0, \ldots, x_{n} ] \subset \Lambda $
that verifies $ \sigma^{ n + \ell(u) + 2 +  \ell(v) - 1 } (z) \in [ y_0, \ldots, y_{m} ] $.
\end{proof}

We establish a sufficient condition to the \hyperlink{dpm}{denseness of periodic measures}.
\begin{proposition}  \label{proposition-sufficient-condition-dpm}
Let $ \left( \Sigma, \sigma \right) $ be a shift over a countable alphabet that satisfies: 
given $ \delta > 0 $ and $ \nu \in \mathcal{M} \left( \Sigma, \sigma \right) $, there exists a finite alphabet subshift $ X_{ \nu, \delta } \subset \Sigma $ fulfilling the specification property such that $ \nu \left( \Sigma \setminus X_{ \nu, \delta } \right) < \delta $.
Then, $ \left( \Sigma, \sigma \right) $ obeys the \hyperlink{dpm}{denseness of periodic measures}.
\end{proposition}

\begin{proof}
Let $ \nu $ be a $ \sigma $-invariant probability measure on $ \Sigma $ and consider the (weak$^{\ast}$) neighborhood
$ V_{\Sigma} \left( \nu; \, f_1, \ldots, f_m; \, \varepsilon \right) $, for some continuous bounded functions $ f_1, \ldots, f_k $ on $ \Sigma $ and for any $ \varepsilon > 0 $.
We must exhibit a probability measure~$ \mu_{\mathrm{per}} $ which is supported on a periodic orbit and verifies, for every $ 1 \leq l \leq m $,
\begin{equation*}
\left| \int_{\Sigma} f_{l} \, d \nu - \int_{\Sigma} f_{l} \, d \mu_{\mathrm{per}} \right| 
<  \varepsilon.
\end{equation*}

Let $ \delta = \nicefrac{ \varepsilon }{ 3 \max_{1 \leq l \leq m} \| f_l \|_\infty } $ and  $ X = X_{ \nu, \delta } $ be the subshift provided by the hypothesis.
From $ \nu \left( \Sigma \setminus X \right) < \nicefrac{\varepsilon}{ 3 \max_{1 \leq l \leq m} \| f_l \|_\infty } $, we obtain
\begin{equation*}
\left| \int_{\Sigma} f_{l} \, d \nu - \int_{X} f_{l} \, d \nu \right|
< \dfrac{\varepsilon}{3}
\end{equation*}
for each $ 1 \leq l \leq m $.
Since $ \left( X, \sigma|_{X} \right) $ is a finite alphabet shift satisfying the specification property and, on the Borel sets of $ X $, the probability $ \mu = \frac{1}{ \nu(X) }  \, \nu $ is a $ \sigma|_{X} $-invariant measure, we may apply Theorem~1 of~\cite{Sig:IM70}, which guarantees the existence of a probability measure $ \mu_{ \mathrm{per} } $ supported on a periodic orbit of $ X $ in the (weak$^{\ast}$) neighborhood
$ V_{X} \left( \mu; \,\nu ( X ) f_1 |_{X}, \ldots, \nu ( X ) f_m |_{X}; \, \nicefrac{\varepsilon}{3} \right) $.
Clearly, 
$ \mu_{ \mathrm{per} } \in \mathcal{M} \left( \Sigma, \sigma \right) $ and
\begin{equation*}
\left| \int_{X} f_{l} \, d \nu - \nu ( X ) \int_{\Sigma} f_{l} \, d \mu_{ \mathrm{per} } \right|
= \left| \int_{X} \nu ( X )  f_{l} \, d \mu - \int_{X}  \nu ( X ) f_{l} \, d \mu_{ \mathrm{per} } \right|
< \dfrac{\varepsilon}{3}
\end{equation*}
for all $ 1 \leq l \leq m $.
It remains to estimate 
\begin{equation*}
\left| \nu ( X ) \int_{\Sigma} f_{l} \, d \mu_{ \mathrm{per} } - \int_{\Sigma} f_{l} \, d \mu_{ \mathrm{per} } \right|
\leq \left| \nu \big( \Sigma \setminus X \big) \right| \, \| f_l \|_\infty 
< \dfrac{\varepsilon}{3},
\end{equation*}
for every $ 1 \leq l \leq m $.
\end{proof}

\begin{lemma}  \label{lemma-non-markov-specification}
For each $ K > 1 $, the finite alphabet subshift $ \Lambda \cap \left( Z_1 \sqcup Z_2 \sqcup \cdots \sqcup Z_{ K } \right)^{ \mathbb{N} } $
satisfies the specification property.
\end{lemma}

\begin{proof}
Given any two words $ u $ and $ w $ of the language of $ \Lambda \cap \left( Z_1 \sqcup Z_2 \sqcup \cdots \sqcup Z_{ K } \right)^{ \mathbb{N} } $, we show that there exists a positive integer $ N = N(K) $ such that these words can be connected by some transitional word $ v $ of length $ N $, meaning that $ u \, v \, w $ belongs to the language of $ \Lambda \cap \left( Z_1 \sqcup Z_2 \sqcup \cdots \sqcup Z_{ K } \right)^{ \mathbb{N} } $.

We can decompose $ u $ and $ w $ in the following forms:
\begin{itemize}
\item
$ u = \dot{u} \, \ddot{u} $ or $ u = \dot{u} \, \ddot{u} \, 1 $, where $ \ddot{u} $ is the subword of $ u $ with largest length whose symbols belong to the same $ Z_{i} $, for $ i > 1 $;
\item
$ w = \ddot{w} \, \dot{w} $ or $ w = 1 \, \ddot{w} \, \dot{w} $, where $ \ddot{w} $ is the subword of $ w $ with largest length whose symbols belong to the same $ Z_{j} $, for $ j > 1 $.
\end{itemize}

Thanks to \ref{item-C4}, each subshift $ \Lambda \cap Z_{k}^{ \mathbb{N} } $ satisfies the specification property, for $ k > 1 $. 
Thus, when $ u = \dot{u} \, \ddot{u} $, there exists a positive integer $ N_{i} $ such that one can find a word~$ \bar{u} $, with $ \ell ( \bar{u} ) = N_{i} $, for which $ \ddot{u} \, \bar{u} \, \max Z_{i} $ belongs to the language of $ \Lambda \cap Z_{i}^{ \mathbb{N} } $.
Similarly, when $ w = \ddot{w} \dot{w} $, there is a word~$ \bar{w} $, with $ \ell( \bar{w} ) = N_{j} $, for which $ \min Z_{j} \, \bar{w} \, \ddot{w} $ belongs to the language of $ \Lambda \cap Z_{j}^{ \mathbb{N} } $.

Therefore, in any situation, we are able to find words $ \check{v} $ and $ \hat{v} $, with lengths in $ \left\{ 0, N_{2} + 2, \ldots, N_{K} + 2 \right\} $, such that $ 1 $ is the last symbol of $ \check{v} $ and the first one of $ \hat{v} $, and both words $ u \, \check{v} $ and $ \hat{v} \, u $ belong to the language of $ \Lambda \cap \left( Z_1 \sqcup Z_2 \sqcup \cdots \sqcup Z_{ K } \right)^{ \mathbb{N} } $. 
Note that words of the form $ \left( 2, \ldots, 2, 3 \right) $, where the symbol $ 2 $ appears $ m \geq 1 $ consecutive times, are admissible in this subshift.
Hence, we can consider $ N (K) = 2 \, \max_{2 \leq i \leq K} N_{i} + 6 $ to conclude the proof.
\end{proof}

\begin{corollary}  \label{corollary-non-markov-dpm}
The set of ergodic probabilities supported on
periodic orbits of $ \Lambda $ is (weak$^{\ast}$) dense among the $ \sigma $-invariant measures.
\end{corollary}

\begin{proof}
Due to \ref{item-C1}, we can decompose the shift space 
$ \Lambda = \Delta \sqcup \, \bigsqcup_{k \geq 0} \Lambda_{k} $, where 
\begin{align*}
\Delta & = \left\{ 
\left( x_0, x_1, \ldots \right) \in \Lambda \,:\,
x_n \neq 1, \text{ for all } n \geq 0 \right\}
\\
\text{and} \qquad
\Lambda_{k} & = \left\{
\left( x_0, x_1, \ldots \right) \in \Lambda\,:\,
x_k = 1 \text{ and } 
x_n \neq 1, \text{ for all } n > k \right\},
\quad \text{for } k \geq 0.
\end{align*}
First, we will prove that $ \nu \left( \bigsqcup_{k \geq 0} \Lambda_{k} \right) = 0 $ for every $ \nu \in \mathcal{M} \left( \Lambda, \sigma \right) $. 
From the Birkhoff Ergodic Theorem, we have
$ \nu \left( \Lambda_k \right) = 
\int \big( \lim_{n \to \infty} \frac{1}{n} \sum_{ i = 0 }^{ n - 1 } \mathds{1}_{ \Lambda_k } \circ \sigma^{i} \big) \, d\nu $, where $ \mathds{1}_{ \Lambda_k } $ is the characteristic function of the set $ \Lambda_k $.
If $ x \in \Delta $, then $ \mathds{1}_{ \Lambda_k } \circ \sigma^{i} (x) = 0 $ for all $ i \geq 0 $.
Moreover, if $ x \in \Lambda_l $ it follows that
\begin{equation*}
\mathds{1}_{ \Lambda_k } \circ \sigma^{i} (x) = \left\{
\begin{array}{cl}
1, & \text{if } i = l - k \\
0, & \text{otherwise}
\end{array}
\right.
\quad \Longrightarrow \qquad
\sum_{ i = 0 }^{ n - 1 } \mathds{1}_{ \Lambda_k } \circ \sigma^{i} (x) = \left\{
\begin{array}{cl}
1, & \text{if } n \geq l - k \\
0, & \text{otherwise}
\end{array}
\right. .
\end{equation*}
Hence, the previous limit is zero everywhere.

It is straightforward from condition \ref{item-C3} that $ \Delta $ is a disjoint union of finite alphabet subshifts $ \Delta_{k} = \Lambda \cap Z_{k}^{ \mathbb{N} } $, for~$ k \geq 2 $. 
In particular, $ \sum_{ k \geq 2 } \nu \left( \Delta_{k} \right) = \nu \left( \Delta \right) = \nu \left( \Lambda \right) = 1 $ for any $ \nu \in \mathcal{M} \left( \Lambda, \sigma \right) $. 

Given $ \delta > 0 $ and a $ \sigma $-invariant probability measure $ \nu $, there exists some integer $ K = K ( \nu,\delta ) $, such that $ \sum_{ k > K } \nu \left( \Delta_{k} \right) < \delta $.
Note that, by \ref{item-C2}, \ref{item-C3} and \ref{item-C4}, the smallest finite alphabet subshift of $ \Lambda $ which contains the subset $ \bigsqcup_{ k = 2 }^{ K } \Delta_{k} $ is given as
\begin{equation*}
X_{ \nu, \delta } = \Lambda \cap \left( Z_1 \sqcup Z_2 \sqcup \cdots \sqcup Z_{ K } \right)^{ \mathbb{N} },
\end{equation*}
which verifies 
$ \nu \left( X_{ \nu, \delta } \right) 
\geq \sum_{ k = 2 }^{K} \nu \left( \Delta_{k} \right)
 \geq 1 - \delta $.
Thanks to Lemma~\ref{lemma-non-markov-specification}, we conclude the result by applying Proposition~\ref{proposition-sufficient-condition-dpm}.
\end{proof}

\begin{Backmatter}

\begin{ack}
E.~Garibaldi, J.~Gomes and M.~Sobottka were partially supported by 
AMSUD240026,
CAPES MathAmSud 88881.985536/2024-01.
\enspace
E.~Garibaldi and J.~Gomes were partially supported by 
ANR-22-CE40-3348 THERMOGAMAS,
FAPESP 24/04685-2.
\enspace
J.~Gomes was partially supported by 
PROEX 88881.973919/2024-01.
\enspace
\end{ack}


\end{Backmatter}

\end{document}

%% file: packs.tex
\usepackage{lmodern}

\theoremstyle{cupplain}

\newtheorem*{theorem-E}{Existence Theorem}

\newtheorem*{theorem-E-SV}{Existence Theorem for Summable Variation Potentials}

\newtheorem*{corollary-LH}{Corollary for Locally H\"older Potentials}

\newtheorem{proposition}[theorem]{Proposition}

\theoremstyle{cupremark}

\usepackage{tikz}
\usetikzlibrary{arrows, automata, backgrounds, chains, matrix, shadows, fit, positioning, calc, scopes, fadings, shapes.geometric, patterns}

\usepackage{enumitem}

\usepackage[nice]{nicefrac}
\usepackage{dsfont}